\documentclass[10pt]{siamltex}
\usepackage{amsfonts,a4wide}
\usepackage{amssymb}
\usepackage{graphicx}
\usepackage{mathtools}
\usepackage{enumerate}
\newtheorem{remark}[theorem]{Remark}
\newtheorem{example}[theorem]{Example}
\numberwithin{equation}{section}
\numberwithin{table}{section}
\numberwithin{figure}{section}
\newcommand{\eproof}{\space
    {\ \vbox{\hrule\hbox{\vrule height1.3ex\hskip0.8ex\vrule}\hrule}}\par}

\newcommand{\IM}{\mbox{\rm Im\,}}
\newcommand{\RE}{\mbox{\rm Re\,}}

\newcommand{\Real}{\mathbb{R}}
\newcommand{\Comp}{\mathbb{C}}
\newcommand{\eps}{\varepsilon}

\newcommand{\To}{\Rightarrow}

\DeclareMathOperator{\ii}{i}

\newcommand{\eq} [1] {\begin{equation}\label{#1}}
\newcommand{\en} {\end{equation}}
\newcommand {\eqn}  {\begin{eqnarray}}
\newcommand {\enn}  {\end{eqnarray}}
\newcommand {\bstar}    {\begin{eqnarray*}}
\newcommand {\estar}    {\end{eqnarray*}}
\newcommand {\mat}  [1] {\left[\begin{array}{#1}}
\newcommand {\rix}      {\end{array}\right]}
\newcommand{\norm}[1]{\left\| #1 \right\|}
\newcommand{\set}[1]{\left\{ #1 \right\}}

\catcode`@=11     
\font\tenex=cmex10 
\newdimen\p@renwd
\setbox0=\hbox{\tenex B} \p@renwd=\wd0 
\def\bmat#1{\begingroup \m@th
  \setbox\z@\vbox{\def\cr{\crcr\noalign{\kern2\p@\global\let\cr\endline}}%
    \ialign{$##$\hfil\kern2\p@\kern\p@renwd&\thinspace\hfil$##$\hfil
      &&\quad\hfil$##$\hfil\crcr
      \omit\strut\hfil\crcr\noalign{\kern-\baselineskip}%
      #1\crcr\omit\strut\cr}}%
  \setbox\tw@\vbox{\unvcopy\z@\global\setbox\@ne\lastbox}%
  \setbox\tw@\hbox{\unhbox\@ne\unskip\global\setbox\@ne\lastbox}%
  \setbox\tw@\hbox{$\kern\wd\@ne\kern-\p@renwd\left[\kern-\wd\@ne
    \global\setbox\@ne\vbox{\box\@ne\kern2\p@}%
    \vcenter{\kern-\ht\@ne\unvbox\z@\kern-\baselineskip}\,\right]$}%
  \null\;\vbox{\kern\ht\@ne\box\tw@}\endgroup}
\catcode`@=12    

\def\diag{\mathop{\mathrm{diag}}}

\newif\ifMDlatex
\catcode`@=11     
\def\MD@us#1{\csname#1style\endcsname}
\def\MD@uf#1{\csname#1font\endcsname}
\def\MD@t{text}
\def\MD@s{script}
\def\MD@ss{scriptscript}
\newdimen\MD@unit
\MD@unit\p@
\def\MD@changestyle#1{
  \relax\MD@unit0.1\fontdimen6\MD@uf{#1}0
  \everymath\expandafter{\the\everymath\MD@us{#1}}
}
\def\MD@dot{$\m@th\ldotp$}
\def\MD@palette#1{\mathchoice{#1\MD@t}{#1\MD@t}{#1\MD@s}{#1\MD@ss}}
\def\MD@ddots#1{{\MD@changestyle{#1}%
  \mkern1mu\raise7\MD@unit\vbox{\kern7\MD@unit\hbox{\MD@dot}}%
  \mkern2mu\raise4\MD@unit\hbox{\MD@dot}%
  \mkern2mu\raise \MD@unit\hbox{\MD@dot}\mkern1mu}}%
\def\MD@iddots#1{{\MD@changestyle{#1}%
  \mkern1mu\raise \MD@unit\hbox{\MD@dot}%
  \mkern2mu\raise4\MD@unit\hbox{\MD@dot}%
  \mkern2mu\raise7\MD@unit\vbox{\kern7\MD@unit\hbox{\MD@dot}}}}%
\def\MD@vdots#1{\vbox{\MD@changestyle{#1}%
    \baselineskip4\MD@unit\lineskiplimit\z@
    \kern6\MD@unit\hbox{\MD@dot}\hbox{\MD@dot}\hbox{\MD@dot}}}%
\ifMDlatex
  \DeclareRobustCommand\ddots{\mathinner{\MD@palette\MD@ddots}}%
  \DeclareRobustCommand\iddots{\mathinner{\MD@palette\MD@iddots}}%
  \DeclareRobustCommand\vdots{\mathinner{\MD@palette\MD@vdots}}%
\else
  \def\ddots{\mathinner{\MD@palette\MD@ddots}}%
  \def\iddots{\mathinner{\MD@palette\MD@iddots}}%
  \def\vdots{\mathinner{\MD@palette\MD@vdots}}%
\fi
\catcode`@=12    

\newcommand {\comment}[1]{} 
\usepackage{color}

\begin{document}

\title{Matrix pencils with coefficients that have positive semidefinite Hermitian part}
\author{C. Mehl \footnotemark[3]~\footnotemark[1]
\and V. Mehrmann\footnotemark[3]~\footnotemark[1]
\and  M. Wojtylak \footnotemark[2]~\footnotemark[4]
}

\maketitle

\begin{abstract}
We analyze when an arbitrary  matrix pencil is
equivalent to a dissipative Hamiltonian pencil and show that this heavily restricts the spectral properties. In order to relax the
spectral properties, we introduce matrix pencils with coefficients that have positive semidefinite Hermitian parts. We will make a
detailed analysis of their
spectral properties and their numerical range. In particular, we relate the Kronecker structure of these pencils to that of an
underlying skew-Hermitian pencil  and discuss their regularity, index, numerical range, and location of eigenvalues.
Further, we study matrix polynomials with positive semidefinite Hermitian coefficients and use linearizations
with positive semidefinite Hermitian parts to derive sufficient conditions for a spectrum in the left half plane and derive bounds
on the index.
%

%
\end{abstract}

{\bf Keywords.} matrix pencils with coefficients that have a positive semidefinite Hermitian part, posH pencils, stability,
pencil regularity,  index structure, dissipative Hamiltonian system,  Kronecker canonical form,
matrix polynomials with positive semidefinite Hermitian coefficients.


\renewcommand{\thefootnote}{\fnsymbol{footnote}}
\footnotetext[3]{
Institut f\"ur Mathematik, MA 4-5, Technische Universit\"at Berlin, Stra{\ss}e des 17. Juni 136,
D-10623 Berlin, Germany.
\texttt{$\{$mehl,mehrmann$\}$@math.tu-berlin.de}.
}

\footnotetext[2]{Instytut Matematyki, Wydzia\l{} Matematyki i Informatyki,
Uniwersytet Jagiello\'nski, Krak\'ow, ul. \L ojasiewicza 6, 30-348 Krak\'ow, Poland
   \texttt{michal.wojtylak@uj.edu.pl}.}
\footnotetext[4]{   Supported by the Alexander von Humboldt Foundation.}
\footnotetext[1]{
Partially supported by {\it Deutsche Forschungsgemeinschaft} through
Priority Program 1984 'Hybride und multimodale Energiesysteme:
 Systemtheoretische Methoden für die Transformation und den Betrieb komplexer Netze'.
}
\renewcommand{\thefootnote}{\arabic{footnote}}

\section{Introduction}
{\color{black}
In this paper we generalize the class of \emph{dissipative Hamiltonian (dH)} matrix pencils, which are pencils  of the form
\begin{equation}\label{dhpencil}
P(\lambda) =\lambda E-A=\lambda E-(J-R)Q,
\end{equation}
where $E,J,R\in\mathbb F^{n,n}$, ($\mathbb F\in\{\mathbb R,\mathbb C\}$), $J^*=-J$, $Q^*E=E^*Q\geq 0$, and $R^*=R\geq 0$.
Here 
${}^*$ stands for the conjugate transpose and 
 $M\geq 0$ (or $M>0$) denotes that the Hermitian matrix $M$ is positive semidefinite (or positive definite, respectively).}
Such dH pencils have many favourable properties, e.g.~all finite eigenvalues are in the closed left half plane and
all purely imaginary eigenvalues are semisimple, except possibly the eigenvalues $0,\infty$ which may have Jordan blocks of size at most
two, see \cite{MehMW18}. For a detailed discussion of dissipative and port-Hamiltonian systems and their applications we refer to
\cite{AltMU21,BeaMV19,BeaMXZ18,GerHH21,MehM19,SchJ14,SchM18}. The observation that the positivity and symmetry structures of the
coefficient matrices of dH pencils lead to these restrictions
in the spectrum shows that the imposed structural conditions - although looking rather simple - are in fact very strong. Imposing
that the matrices $E,A$ are real leads to several further spectral properties, see \cite{MehMW21}.

Since spectral properties are invariant under equivalence transformations of the matrix pencil, but the dissipative Hamiltonian
structure is not, it is clear that there are many pencils that have the same spectral properties but that do not have the
structure of the pencil as in~\eqref{dhpencil}. As our first result (Theorem~\ref{stable=dh}) we will characterize  when a general
matrix pencil
$L(\lambda)=\lambda E-A$ is equivalent to a dH pencil and we will show that it is
necessary and sufficient that the mentioned spectral properties hold.

{In several applications matrix pencils arise that carry a structure that is related to, but more general than the one of dH pencils.
These are square pencils of the form
%
\begin{equation}\label{jj-pencil}
\lambda (J_1+R_1)+(J_2+R_2),\quad \text{ with }J_i=-J_i^*,\ R_i\geq 0,\ i=1,2.
\end{equation}
In other words, we will assume that the (uniquely defined)  Hermitian part $R_i$ of each coefficient  is positive semidefinite.
We will call these  pencils \emph{posH pencils}, abbreviating  `positive semidefinite Hermitian part coefficients'. }
If $J_1=0$ (or $J_2=0$) then the pencil in~\eqref{jj-pencil} simply reduces to a dH pencil as in~\eqref{dhpencil} with $Q=I$
(or its reversal, respectively) and thus, all eigenvalues of the pencil are in the closed left half plane.
This is no longer true if both $J_1$ and $J_2$ are nonzero as the following example shows.

\begin{example}\label{ex:unstable}{\rm
Consider the pencil
\[
\lambda \begin{bmatrix}  0 & -1 & 0 \\ 1 & 0 & 0 \\ 0 & 0 & 1 \end{bmatrix}+
\begin{bmatrix} 1 & 0 & 0 \\ 0 & 0 & 1 \\ 0 & -1 & 0 \end{bmatrix}
\]
which has the form~\eqref{jj-pencil} and is a linearization of the scalar polynomial $P(\lambda)=\lambda^3 +1$ that has two roots with
positive real part.
}
\end{example}

In view of Example~\ref{ex:unstable} one may initially think that the structure of posH pencils is rather weak compared to that of
dH pencils, but we will show in this paper that posH pencils still have many special properties. Furthermore, they are of great
importance in applications which makes it necessary to analyze and study them in detail. Let us give two motivating examples.

\begin{example}\label{ex:MooreGibbsThomson}{\rm\color{black}
The space  discretization of the Moore-Gibbs-Thompson equation \cite{Ben21,KalN19}  leads to cubic matrix polynomials
\[
P(\lambda) =\sum_{i=0}^3 \lambda^iA_i,
\]
where all coefficients $A_i$ are real symmetric and positive definite. Using structured linearization (see Theorem~\ref{Benner-explained}
for details) 
one obtains a pencil of the form~\eqref{jj-pencil} given by
\begin{equation}\label{skewlin}
\lambda \begin{bmatrix}  0 & -A_3 & 0 \\ A_3 & A_2 & 0 \\ 0 & 0 & A_0 \end{bmatrix}+
\begin{bmatrix} A_3 & 0 & 0 \\ 0 & A_1 & A_0 \\ 0 & -A_0 & 0 \end{bmatrix}.
\end{equation}
We will analyze under
which conditions all eigenvalues of this pencil are in the open left half-plane, see Theorem~\ref{Benner-explained} and
Corollary~\ref{Benner-explained2}.}
\end{example}

Example~\ref{ex:MooreGibbsThomson} illustrates that posH pencils arise as linearizations of higher order matrix polynomials
with positive (semi-)definite Hermitian coefficients. {\color{black} Further constructions of this type can be found in
Remarks~\ref{rem:5.4.21} and \ref{rem:5.4.21.2}.} However, there are other situations that can be modeled with the help of
posH pencils.

\begin{example}\label{ex:brake}{\rm  In the analysis of disk brake squeal, see \cite{GraMQSW16}, one has to analyze the spectral
properties of quadratic matrix polynomials $\lambda^2 M + \lambda(D+G) +K+N$, with real symmetric positive semidefinite matrices
$M, D, K$ and real skew-symmetric matrices $G,N$. Brake squeal is associated to a flutter instability arising at the brake-pad disk
interface and it is correlated to eigenvalues with positive real part.
Consider the linearization 
\[
\lambda
\begin{bmatrix}
M & 0 \\ 0 & K- N
\end{bmatrix}  + \begin{bmatrix} D+G & K+N \\ -K+N & 0 \end{bmatrix}
\]
that has the form~\eqref{jj-pencil}. If the contribution from the skew-symmetric matrix  $N$ is zero then this is a dH pencil and if the
norm of $N$ is sufficiently small, then this pencil still has all eigenvalues in the left half plane.
However, if the norm of $N$ is larger, then eigenvalues in the right half complex plane occur that may lead to brake squeal.
}
\end{example}

The three presented examples show that extra assumptions for pencils of the form~\eqref{jj-pencil} are needed to guarantee that all
eigenvalues of such pencils or related matrix polynomials are in the left half plane. We will derive such conditions and also
analyze general spectral properties.

The paper is organized as follows.
In Section~\ref{sec:prelim} we present some preliminary results and introduce relevant notation. In Section~\ref{sec2a} we discuss
necessary and sufficient conditions for a pencil to be equivalent to a dH pencil as in~\eqref{dhpencil}.

In Section~\ref{sec3} we relate the Kronecker structure of a posH pencil of the form $\lambda (J_1+R_1)+J_2+R_2$ to that
of the underlying skew-Hermitian pencil $\lambda J_1+J_2$ with particular emphasis on regularity, the index of the pencil, and
positive eigenvalues. {The key result here is Theorem~\ref{thm:chains}, which says that the singular part of a posH pencil
in  \eqref{jj-pencil} is contained in the common kernel of $R_1$ and $R_2$ and in the singular part of $\lambda J_1+J_2$.
This fact leads to several necessary or sufficient conditions for regularity of posH pencils, see Corollaries~\ref{cor:reg1}
and \ref{cor2}, and Theorem~\ref{rem:nocommon}.

In Subsection~\ref{sec:commoniso} we first consider the numerical range for posH pencils, in particular we link the existence
of common isotropic vectors with regularity of the pencil, see Theorem~\ref{rem:nocommon}. In Subsection~\ref{sec:pack} we
localize the numerical range and the spectrum in a pacman-like shape, the main result is Theorem~\ref{3/4}. In
Subsection~\ref{sec:left} we provide several sufficient conditions that guarantee that the numerical range
or at least the spectrum of a posH pencil is contained in the closed left half plane - a condition that is necessary
for stability of the pencil.

In Section~\ref{sec5} we  consider the special case of matrix polynomials with positive semidefinite Hermitian
coefficients, i.e., the skew-Hermitian parts of the coefficients are all zero.
We analyze their index in Theorem~\ref{thm:indexd} and localize the spectrum in Theorem~\ref{Benner-explained}. This is done
by showing that these polynomials can be linearized by posH matrix pencils and by applying the results from previous sections.
}

\section{Preliminaries}\label{sec:prelim}
We denote by $\mathbb F^{n,m}[\lambda]$ the set of matrix polynomials with coefficients in the set $\mathbb F^{n,m}$ of $n\times m$
matrices over $\mathbb F$ ($\mathbb F\in\{\mathbb R,\mathbb C\}$).
For a pencil $L(\lambda) =\lambda E-A\in  \mathbb F^{n,m}[\lambda]$ the \emph{reversal} $\operatorname{rev}L(\lambda)$ is defined as
$\operatorname{rev}L(\lambda)=\lambda A-E$. Two pencils $L(\lambda), \tilde L(\lambda) \in \mathbb F^{n,m}[\lambda]$ are called equivalent
if there exists invertible matrices $S\in \mathbb F^{n,n}$, $T\in \mathbb F^{m,m}$ such that $L(\lambda)=S \tilde L(\lambda) T$.
To analyze the spectral properties of matrix pencils we will employ the \emph{Kronecker canonical form} \cite{Gan59a}. Denote by
$\mathcal J_k(\lambda_0)$ the standard \emph{upper triangular Jordan block} of size $k\times k$ associated with the eigenvalue
$\lambda_0$ and let $\mathcal L_k$ denote the standard \emph{right Kronecker block} of size $k\times(k+1)$, i.e.,
\[
\mathcal L_k=\lambda\left[\begin{array}{cccc}
1&0\\&\ddots&\ddots\\&&1&0
\end{array}\right]-\left[\begin{array}{cccc}
0&1\\&\ddots&\ddots\\&&0&1
\end{array}\right]\quad\mbox{and}\quad \mathcal J_k(\lambda_0)=\left[\begin{array}{cccc}
\lambda_0&1\\&\ddots&\ddots\\&&\ddots&1\\&&&\lambda_0\end{array}\right].
\]

\begin{theorem}[Kronecker canonical form]\label{th:kcf}
Let $E,A\in {\mathbb C}^{n,m}$. Then there exist nonsingular matrices
$S\in {\mathbb C}^{n,n}$ and $T\in {\mathbb C}^{m,m}$ such that
\begin{equation}\label{kcf}
S(\lambda E-A)T=\diag({\cal L}_{\epsilon_1},\ldots,{\cal L}_{\epsilon_p},
{\cal L}^\top_{\eta_1},\ldots,{\cal L}^\top_{\eta_q},
{\cal J}_{\rho_1}^{\lambda_1},\ldots,{\cal J}_{\rho_r}^{\lambda_r},{\cal N}_{\sigma_1},\ldots,
{\cal N}_{\sigma_s}),
\end{equation}
where the parameters $p,q,r,s,\epsilon_1,\dots,\epsilon_p,\eta_1,\dots,\eta_q,\rho_1,\dots,\rho_r,\sigma_1,\dots,\sigma_s$
are nonnegative integers,
$\lambda_1,\dots,\lambda_r\in\mathbb C$, and ${\cal J}_{\rho_i}^{\lambda_i}=I_{\rho_i}-\mathcal J_{\rho_i}(\lambda_i)$
for $i=1,\dots,r$ as well as $\mathcal N_{\sigma_j}=\mathcal J_{\sigma_j}(0)-I_{\sigma_j}$ for $j=1,\dots,s$.
This form is unique up to permutation of the blocks.
\end{theorem}

For real matrices a real version of the Kronecker canonical form is obtained under real transformation matrices $S,T$.
In this case the  blocks ${\cal J}_{\rho_j}^{\lambda_j}$ with $\lambda_j=\alpha_j+\ii\beta_j\in\Comp\setminus\Real$ have to be replaced
with corresponding blocks in \emph{real Jordan canonical form} with diagonal blocks of the form
\[
J_{\rho_j}(\alpha,\beta):=\left[\begin{array}{cccc}
\Lambda_j&I_2\\&\ddots&\ddots\\&&\ddots&I_2\\&&&\Lambda_j\end{array}\right]\in\mathbb R^{2\rho_j,2\rho_j},
\quad \Lambda_j:=\mat{cc}\alpha_j&\beta_j\\ -\beta_j&\alpha_j\rix
\]
associated to the corresponding pair of conjugate complex eigenvalues $\alpha_j\pm \ii\beta_j$, but the other
blocks have the same structure as in the complex case.

An eigenvalue is called \emph{semisimple} if the largest associated Jordan block has size one.
The sizes $\eta_j$ and $\epsilon_i$ of the rectangular blocks
are called the \emph{left and right minimal indices} of $\lambda E-A$, respectively.
If $\eta\geq 0$ is a left minimal index, then there exists a singular chain $(x_1,\dots,x_{\eta+1})$ of
vectors satisfying $x_1^\top E=0$, $x_{j+1}^\top E=x_j^\top A$, $j=1,\dots,\eta$ and $x_{\eta+1}A=0$.
Similarly, if $\epsilon\geq 0$ is a right minimal index, then there exists a singular chain $(x_1,\dots,x_{\epsilon+1})$ of
vectors satisfying $Ex_1=0$, $Ex_{j+1}=Ax_j$, $j=1,\dots,\epsilon$ and $Ax_{\epsilon+1}=0$.
The matrix pencil $\lambda E-A$, $E,A \in \mathbb C^{n,m}$ is called \emph{regular} if $n=m$ and
$\operatorname{det}(\lambda_0 E-A)\neq 0$ for some $\lambda_0 \in \mathbb C$,
otherwise it is called \emph{singular}.
A pencil is singular if and only if it has blocks of at least one of the types ${\cal L}_{\eps_j}$ or
${\cal L}^\top_{\eta_j}$ in the Kronecker canonical form.

The values $\lambda_1,\dots,\lambda_r\in\mathbb C$ are called the finite eigenvalues of $\lambda E-A$. If $s>0$, then
$\lambda_0=\infty$ is said to be an eigenvalue of $\lambda E-A$. (Equivalently, zero is then an eigenvalue of
the reversal $\lambda A-E$ of the pencil $\lambda E-A$.)

The sum of all sizes of blocks that are associated with a fixed eigenvalue
$\lambda_0\in\mathbb C\cup\{\infty\}$ is called the \emph{algebraic multiplicity} of $\lambda_0$, while the individual
sizes of the Jordan blocks are called the \emph{partial multiplicities} of $\lambda_0$.
The size of the largest block ${\cal N}_{\sigma_j}$ is
called the \emph{index} $\nu$ of the pencil $\lambda E-A$, where, by convention,  $\nu=0$ if $E$ is invertible.

The pencil is called {\em stable} if it is regular, if all eigenvalues are in the closed left half plane, and if the ones lying on the
imaginary axis (including infinity) are semisimple. Otherwise the pencil is called {\em unstable}.

\section{Pencils that are equivalent to dH pencils}\label{sec2a}
It is a natural question to ask under which conditions a posH pencil is equivalent to a dH pencil. It turns out that
the answer is obtained by a general characterization including matrix pencils without special symmetry and positivity
structures. Parts of the following
result  were discovered independently in \cite{FauMPSW21}.
\begin{theorem}\label{stable=dh}
\begin{enumerate}[\rm (i)]
\item\label{Qarb} A pencil $L(\lambda)\in \mathbb F^{n,n}[\lambda]$ is equivalent to a pencil of the form
$P(\lambda)=\lambda E-(J-R)Q$ as in~\eqref{dhpencil} with $\lambda E-Q$ being regular
if and only if the following conditions are satisfied:
\begin{enumerate}[\rm (a)]
\item The spectrum of $L(\lambda)$ is contained in the closed left half plane.
\item The finite nonzero eigenvalues on the imaginary axis are semisimple and the partial multiplicities of the eigenvalue zero are at most two.
\item The index of $L(\lambda)$ is at most two.
\item The left minimal indices are all zero and the right minimal indices are at most one (if there are any).
\end{enumerate}
\item\label{QI} A pencil $L(\lambda)\in \mathbb F^{n,n}[\lambda]$ is equivalent to a pencil of the form $P(\lambda)=\lambda E-(J-R)$
as in~\eqref{dhpencil} (i.e., with $Q=I$) if and only if the following conditions are satisfied:
\begin{enumerate}[\rm (a)]
\item The spectrum of $L(\lambda)$ is contained in the closed left half plane.
\item The finite eigenvalues on the imaginary axis (including zero) are semisimple.
\item The index of $L(\lambda)$ is at most two.
\item The left and right minimal indices are all zero (if there are any).
\end{enumerate}
\end{enumerate}
\end{theorem}
\proof
The ``only if'' direction for \eqref{Qarb} was proved in \cite{MehMW18} and the one for (ii) in \cite{MehMW21},
see also \cite{GilMS18} for the matrix case.

For the ``if'' direction, we may assume
without loss of generality that $L(\lambda)$ is in Kronecker canonical form.
In particular, we may consider each block separately. First, we prove \eqref{QI} and we distinguish the following cases for
$\lambda_0= \alpha + \ii\beta$ with $\alpha,\beta \in \mathbb R$.

{\it Case 1)}: $\mathbb F=\mathbb C$.\\
\qquad{\it Subcase 1a)}: $L(\lambda)=\lambda I_n-J_n(\lambda_0)$, $n\geq 1$, with $\alpha<0$.
 Then $L(\lambda)$ is equivalent to the pencil $\lambda I_n-M=\lambda E-(J-R)Q$, with
\[
M=\mat{cccc}\alpha+\ii\beta&\alpha&&0\\ &\alpha+\ii\beta&\ddots&\\ &&\ddots&\alpha\\ 0&&&\alpha+\ii\beta\rix,\ E=Q=I_n,
\]
\[
J=\frac{1}{2}(M-M^*)=\mat{cccc}\!\!\!\ii\beta&\frac{1}{2}\alpha&&\\ \!\!\!-\frac{1}{2}\alpha&\ii\beta&\ddots&\\
&\ddots&\ddots&\frac{1}{2}\alpha\!\\ &&-\frac{1}{2}\alpha&\ii\beta\!\rix,\;\,
R=-\frac{1}{2}(M+M^*)=-\mat{cccc}\!\alpha&\frac{1}{2}\alpha&&\\ \!\frac{1}{2}\alpha&\alpha&\ddots&\\
&\ddots&\ddots&\frac{1}{2}\alpha\!\\ &&\frac{1}{2}\alpha&\alpha\!\rix.
\]
By \cite[Proposition 2.2]{AndD11} which is a combination of Theorem 2.4 and Proposition 2.5 in \cite{JohNT96},
it follows that $R\geq 0$ (in fact $R>0)$.\\
\quad {\it Subcase 1b)}: $L(\lambda)=\lambda-\ii\beta$ with $\beta\in\mathbb R$. Here, we have $E=Q=1$, $J=\ii\beta$ and $R=0$.\\
\quad{\it Subcase 1c)}: $L(\lambda)=\lambda 0-1$. This pencil is equivalent to $\lambda 0+1$ and we can take
$E=J=0$, and $R=Q=1$.\\
\quad {\it Subcase 1c')}: $L(\lambda)=\lambda J_2(0)-I_2$. Then $L(\lambda)$ is equivalent to the pencil
\[
\lambda E-JQ\quad\mbox{with}\quad E=\mat{cc}0&0\\ 0&1\rix,\quad J=\mat{cc}0&1\\ -1&0\rix,\quad R=0\quad\mbox{and}\quad Q=I_2.
\]
\quad {\it Subcase 1d)}: Since the pencil is square the numbers of left and right minimal indices are equal and hence
each pair corresponds in the Kronecker canonical form to a block $L(\lambda)=\lambda 0-0$. Here we can take $Q=1$ and $E=R=J=0$.

{\it Case 2)}: $\mathbb F=\mathbb R$.\\
\quad{\it Subcase 2a)}: $L(\lambda)=\lambda I_n-J_n(\alpha)$ with $\alpha<0$. This case works exactly as
Subcase 1a) with $\beta=0$. \\
\quad {\it Subcase 2a')}: $L(\lambda)=\lambda I_{2n}-J_n(\alpha,\beta)$ with $\alpha<0$ and $\beta\neq 0$. In this case $L(\lambda)$
is equivalent to the matrix pencil $\lambda I_{2n}-M=\lambda E-(J-R)Q$ with
\begin{eqnarray*}
M&=&\left[\begin{array}{cccc}
\Lambda&\alpha I_2\\&\ddots&\ddots\\&&\ddots&\alpha I_2\\&&&\Lambda\end{array}\right],
\quad \Lambda=\mat{cc}\alpha&\beta\\ -\beta&\alpha\rix,\\
E=Q=I_{2n},&& J=\frac{1}{2}(M-M^T),\quad R=-\frac{1}{2}(M+M^T).
\end{eqnarray*}
 Again, by combining Theorem 2.4 and Proposition 2.5 in \cite{JohNT96} it follows that
$R\geq 0$.\\
\quad {\it Subcase 2b)}: $L(\lambda)=\lambda I_{2}-J_1(0,\beta)$ with $\beta\neq 0$. Here we have $L(\lambda)=\lambda E-(J-R)Q$ with
$E=Q=I_2$, $R=0$, and $J=J_1(0,\beta)$.\\
\quad {\it Subcase 2b')}: $L(\lambda)=\lambda -0$. Here we have $E=Q=1$ and $J=R=0$.\\
\quad The subcases 2c) and 2d) are identical to the subcases 1c) and 1d) as the corresponding matrices are all real.

To prove \eqref{Qarb} it remains to consider one additional block of the form
$L(\lambda)=\lambda I_2-J_2(0)$ in subcase 1b) and one additional combination of minimal indices in subcase 1d).
In the first case we have $L(\lambda)=\lambda E-(J-R)Q$ with
\[
E=I_2,\ R=0,\ J=\mat{cc}0&1\\ -1&0\rix,\  Q=\mat{cc}0&0\\ 0&1\rix.
\]
As all matrices are real, this case works for both the real and the complex case.

In the second case, note that again the numbers of left and right minimal indices must be equal. For a pair of left and right
minimal indices $(0,0)$ we are in subcase 1d) of~\eqref{Qarb}. For a pair of left and right minimal indices $(0,1)$ we have
a block
\[
L(\lambda)=\lambda\mat{cc}1&0\\ 0&0\rix-\mat{cc}0&1\\ 0&0\rix
\]
in the Kronecker canonical form. Here, we can take
\[
E=\mat{cc}1&0\\ 0&0\rix,\;R=0,\; J=\mat{cc}0&1\\ -1&0\rix\;\quad\mbox{and}\quad Q=\mat{cc}0&0\\ 0&1\rix.
\]
Again all matrices are real, so this case works for both the real and the complex case.
\eproof
Theorem~\ref{stable=dh} clearly shows that the spectral properties are precisely characterizing the equivalence to matrix pencils
of the form~\eqref{dhpencil}, so we cannot expect similarly nice spectral properties if we generalize to pencils of
the form \eqref{jj-pencil}.  However, we still get spectral restrictions for such pencils, some of which are associated with
the \emph{numerical range} which is an important tool in investigating stability of matrices, matrix pencils or matrix polynomials.
These will be discussed in the following sections.
\section{On the Kronecker structure of posH matrix pencils}\label{sec3}
In this section we will investigate the Kronecker structure of posH pencils, i.e., matrix pencils of the form \eqref{jj-pencil}.
We have already seen  in Example~\ref{ex:unstable} that posH pencils may have eigenvalues in the right half plane including eigenvalues
on the positive real axis. In fact, without posing further restrictions on the pencil, any eigenvalue in the complex
plane is possible.
\begin{example}\label{ex:jja}\rm
Let $\alpha,\beta\in\mathbb R$. If $\beta\geq 0$ then
$P_1(\lambda)=\ii\lambda+(\beta-\ii\alpha)$, i.e., $J_1=\ii$, $R_1=0$, $J_2=-\ii\alpha$, $R_2=\beta$ is a complex posH matrix pencil
having the eigenvalue $\alpha+\ii\beta$. (If $\beta<0$ then consider the complex posH pencil
$P_2(\lambda)=-P_1(\lambda)$ instead.) In particular, if $\beta=0$ and $\alpha>0$, then $P_1(\lambda)$ is an example of a posH
pencil with an eigenvalue on the positive real axis.
\end{example}

\begin{example}\label{ex:jjb}\rm
For an example with real matrix coefficients consider the posH matrix pencil
$P_r(\lambda)=\lambda (J_1+R_1)+(J_2+R_2)$ with
\[
R_1=0,\quad J_1=\mat{cc}0&1\\ -1&0\rix,\quad J_2=\mat{cc}0&\alpha\\ -\alpha&0\rix
\quad\mbox{and}\quad R_2=\mat{cc}\beta &0\\ 0&\beta\rix,
\]
where $\beta\geq 0$ and $\alpha\in\mathbb R$.
Then $P_r(\lambda)$ has a pair of conjugate complex eigenvalues $\alpha\pm \ii\beta$. In particular, if  $\beta=0$ and $\alpha>0$,
then $P_r(\lambda)$ has a double eigenvalue on the positive real axis.
\end{example}

Although the spectrum may contain any value of the complex plane, not any Kronecker structure is possible for posH pencils.
In the following we will discuss restrictions on the index and the structure of the singular part of such pencils.
We start with two technical results on values $\lambda_0$ and vectors $x$ satisfying
$P(\lambda_0)x=0$. Note that the pencil is not excluded to be singular, so $\lambda_0$ is not necessarily an eigenvalue of $P(\lambda)$.

\begin{lemma}\label{lem:prop2}
Let $P(\lambda)=\lambda (J_1+R_1)+(J_2+R_2)\in\mathbb C^{n,n}[\lambda]$ be a pencil as in~\eqref{jj-pencil} and let $x\in\mathbb C^n$
and $\lambda_0\in\mathbb C$.
\begin{enumerate}[\rm (i)]
\item\label{c-i} If $\,\IM\lambda_0=0$, $\RE\lambda_0>0$ and $P(\lambda_0)x=0$, then
$$
R_1x=R_2x=0\quad\mbox{and}\quad (\lambda_0 J_1+J_2)x=0.
$$
\item\label{c-iii} If $\lambda_0\in\Comp$, $x^*(J_1+R_1)x=0$ and $x^*P(\lambda_0)x=0$ then
 $$
R_1x=R_2x=0\quad \text{ and}\quad  x^*J_1x=0=x^*J_2x.
$$
\end{enumerate}
\end{lemma}
\proof
Let $\lambda_0=\alpha+\ii\beta$ with $\alpha,\beta\in\mathbb R$. First observe that $x^*P(\lambda_0)x=0$ in both cases \eqref{c-i}
and \eqref{c-iii}. Taking the real and imaginary parts independently yields the equations
\begin{eqnarray}
\alpha x^*R_1x+\ii\beta  x^*J_1x&=&-x^*R_2x,\label{16.3.21.1}\\
\alpha x^*J_1x+\ii\beta  x^*R_1x&=&-x^*J_2x\label{16.3.21.2}
\end{eqnarray}

\eqref{c-i} Assume   $\alpha>0$ and $\beta=0$, then we obtain from~\eqref{16.3.21.1}
that $\alpha x^*R_1x=-x^*R_2x$ which, by the semidefiniteness of $R_1$ and $R_2$, is only possible if $R_1x=0=R_2x$. But then
we have $0=P(\lambda_0)x=\lambda_0J_1x+J_2x$.

\eqref{c-iii} Let $\alpha,\beta\in\mathbb R$  be arbitrary.
Then due to $x^*(J_1+R_1)x=0$ one has $x^*J_1x=0$ and $x^*R_1 x=0$, which implies  $R_1x=0$. Hence, thanks to~\eqref{16.3.21.1}
we have $R_2x=0$ and furthermore $x^*J_2x=0=x^*J_1x$ by~\eqref{16.3.21.2}.
\eproof
By Theorem~\ref{stable=dh} the left and right minimal indices of a singular dH pencil  with $Q=I$ as in~\eqref{dhpencil} can
only be zero. This is no longer true for posH pencils of the form~\eqref{jj-pencil}, but the following result
shows that the singular part of posH pencils is still restricted.
\begin{theorem}\label{thm:chains}
Let $P(\lambda)=\lambda (J_1+R_1)+(J_2+R_2)\in\mathbb F^{n,n}[\lambda]$ be a pencil of the form~\eqref{jj-pencil}.
If $(x_1,\dots,x_{k+1})$ is a singular chain associated with
a left or right minimal index $\eta=k$, then $x_1,\dots,x_{k+1}\in\ker R_1\cap\ker R_2$ and $(x_1,\dots,x_{k+1})$ is also
a singular chain of $\lambda J_1+J_2$ associated with a left respectively right minimal index $\eta=k$. 
%
\end{theorem}
\proof
Let $(x_1,\dots,x_{k+1})$ be a singular chain associated with a left or right minimal index $\eta=k$ of $P(\lambda)$. Without
loss of generality, let this be a right minimal index, otherwise, consider the pencil with coefficients that are the conjugate
transposes of that of $P(\lambda)$. Then we have
\begin{equation}\label{chain}
(J_1+R_1)x_1=0,\quad (J_1+R_1)x_{i+1}=(J_2+R_2)x_i,\;i=1,\dots,k,\quad (J_2+R_2)x_{k+1}=0,
\end{equation}
or equivalently, using that $J_i^*=-J_i$ and $R_i^*=R_i$ and multiplying by $-1$,
\begin{equation}\label{chain2}
x_1^*(J_1-R_1)=0,\quad x_{i+1}^*(J_1-R_1)=x_i^*(J_2-R_2),\;i=1,\dots,k,\quad x_{k+1}^*(J_2-R_2)=0.
\end{equation}
We first prove by induction that $R_1x_j=0=R_2x_j$ for all $j=1,\dots,k+1$. From $x_1^*(R_1-J_1)x_1=0$ we get $R_1x_1=0$ and thus
$J_1x_1=(J_1+R_1)x_1=0$. If $k=0$, then we have $x_{k+1}=x_1$ and also $R_2x_1=0$ follows similarly from $(J_2+R_2)x_1=0$.
Otherwise we have
\[
x_1^*(J_2+R_2)x_1=x_1^*(J_1+R_1)x_{2}=0
\]
which implies that $R_2x_1=0$.

Suppose that for some $\ell\geq 1$ we have shown $R_1x_j=0=R_2x_j$
for all $j=1,\dots,\ell$. If $\ell=k$ then we are done, because similar to the previous  argument we then get
$J_2x_{k+1}=0=R_2x_{k+1}$ from $(J_2+R_2)x_{k+1}=0$ and $R_1x_{k+1}=0$ from $x_{k+1}^*(J_1+R_1)x_{k+1}=x_{k+1}^*(J_2+R_2)x_k=0$.

Hence, we may assume that $\ell<k$ and thus $\ell+2\leq k+1$. Using~\eqref{chain} and~\eqref{chain2} we obtain
\[
x_{\ell+1}^*(J_1+R_1)x_{\ell+1}=x_{\ell+1}^*(J_2+R_2)x_\ell=x_{\ell+1}^*(J_2-R_2)x_\ell=x_{\ell+2}^*(J_1-R_1)x_\ell
=x_{\ell+2}^*(J_1+R_1)x_\ell,
\]
where we have used that $R_2x_\ell=0$ and $R_1x_\ell=0$. {We repeat this procedure  $m$ times, obtaining
$x_{\ell+1}^*(J_1+R_1)x_{\ell+1}=x_{\ell+1+m}^*(J_1+R_1)x_{\ell+1-m}$, and we may proceed until $m=\min\set{\ell,k-\ell}$.
If $m=\ell$, i.e., if $\ell+1-m=1$, then  we have
\[
x_{\ell+1}^*(J_1+R_1)x_{\ell+1}=x_{\ell+1+m}^*(J_1+R_1)x_{1}=0,
\]
 while if $m=k-\ell<\ell$, then  we get
\[
x_{\ell+1}^*(J_1+R_1)x_{\ell+1}=x_{k+1}^*(J_1+R_1)x_{\ell+1-m}=x_{k+1}^*(J_2+R_2)x_{\ell-m}=0.
\]}
Thus, in both cases, we finally obtain $x_{\ell+1}^*(J_1+R_1)x_{\ell+1}=0$ which implies that $R_1x_{\ell+1}=0$.

On the other hand, using~\eqref{chain},~\eqref{chain2} and that we just proved that $R_2x_\ell=0$ and
$R_1x_{\ell+1}=0$ hold, we obtain that
\begin{eqnarray*}
x_{\ell+1}^*(J_2+R_2)x_{\ell+1}&=&x_{\ell+1}^*(J_1+R_1)x_{\ell+2}=x_{\ell+1}^*(J_1-R_1)x_{\ell+2}\\
&=&x_{\ell}^*(J_2-R_2)x_{\ell+2}=x_{\ell}^*(J_2+R_2)x_{\ell+2}.
\end{eqnarray*}
As before, {\color{black} we repeat this step $m$ times, obtaining
$x_{\ell+1}^*(J_2+R_2)x_{\ell+1}=x_{\ell+1-m}^*(J_2+R_2)x_{\ell+1+m}$, and we proceed until $m=\min\set{k-\ell,\ell}$. If $m=k-\ell$, then
\[
x_{\ell+1}^*(J_2+R_2)x_{\ell+1}=x_{\ell+1-m}^*(J_2+R_2)x_{k+1}=0,
\]
otherwise we have  $m=\ell<k-\ell$, which gives
\[
x_{\ell+1}^*(J_2+R_2)x_{\ell+1}=x_{1}^*(J_2+R_2)x_{\ell+1+m}=x_{1}^*(J_1+R_1)x_{\ell+2+m}=0.
\]
In both cases, we obtain $x_{\ell+1}^*(J_2+R_2)x_{\ell+1}=0$ which implies that $R_2x_{\ell+1}=0$.

Thus, using an induction argument, we obtain $R_1x_j=0=R_2x_j$ for all $j=1,\dots,k+1$. Inserting that into~\eqref{chain}, we get
\[
J_1x_1=0,\quad J_1x_{i+1}=J_2x_i,\;i=1,\dots,k,\quad J_2x_{k+1}=0
\]
which shows that $(x_1,\dots,x_k)$ is a singular chain of the pencil $\lambda J_1+J_2$ associated with the right minimal
index $\eta=k$.
\eproof

Since a pencil of skew-Hermitian matrices has equal left and right minimal indices, see~\cite{Tho91}, we immediately
obtain by Theorem~\ref{thm:chains} that the same is true for posH pencils.

\begin{corollary}\label{cor:reg1}
Let $P(\lambda)=\lambda (J_1+R_1)+(J_2+R_2)\in \mathbb F^{n,n}[\lambda]$ be a pencil of the form~\eqref{jj-pencil}.
Then the ordered lists of left and right minimal indices of $P(\lambda)$ coincide.
\end{corollary}

In \cite{MehMW21} it was shown that a dH pencil of the form~\eqref{dhpencil} (with $Q=I$) is singular if and only if the three matrices
$E$, $J$, and $R$ have a common kernel. A corresponding result for posH pencils is only true under additional assumptions.
\begin{corollary}\label{thm3}
Let $P(\lambda)=\lambda (J_1+R_1)+(J_2+R_2)\in\mathbb F^{n,n}[\lambda]$ be a pencil of the form~\eqref{jj-pencil}.
If $\ker(J_1)\cap\ker(J_2)\neq\set0$ and if all minimal indices of $\lambda J_1+J_2$ are zero, then
$P(\lambda)$ is singular if and only if the four matrices $J_1,J_2,R_1,R_2$ have a common kernel. Moreover, in this case all left and
right minimal indices of $P(\lambda)$ are zero.
\end{corollary}
\proof
This is a direct consequence of Theorem~\ref{thm:chains}. \eproof

In fact, Corollary~\ref{thm3} is a direct generalization of the corresponding result on dH pencils.
Indeed, if $J_1=0$, then the pencil $\lambda J_1+J_2$ can only have left and right minimal indices equal to zero
and hence, the same is true for any pencil of the form $\lambda R_1+(J_2+R_2)$ with Hermitian positive semidefinite
$R_1$ and $R_2$, i.e., a pencil of the form as in~\eqref{dhpencil} with $Q=I$, see part (v) of \cite[Theorem 2]{MehMW21}.

The latter result on dH pencils can also be generalized to posH pencils in a different way by considering other combinations of three of the four coefficients. Furthermore, by considering pencils built of two of the four coefficients of posH pencils, one can characterize situations when pencils of the form \eqref{jj-pencil} may or may not have positive real eigenvalues.
\begin{corollary}\label{cor2}
Let $P(\lambda)=\lambda (J_1+R_1)+(J_2+R_2)\in\mathbb F^{n,n}[\lambda]$ be a pencil of the form~\eqref{jj-pencil}.
\begin{enumerate}[\rm (i)]
\item\label{sing5} {\it If $P(\lambda)$ is singular then the matrices in each triple $(J_i,R_1,R_2)$, $i=1,2$, have a common kernel.}
\item\label{sing0} {\it If the pencil $\lambda R_1+R_2$ is regular, then $P(\lambda)$ is regular and has no eigenvalues on the real
positive axis.}
\item\label{singii} {\it If the pencil $\lambda J_1+J_2$ is regular, then $P(\lambda)$ is regular and
every real positive eigenvalue of $P(\lambda)$ is also an eigenvalue of the pencil $\lambda J_1+J_2$.}
\item\label{sing3} {\it If the pencil $\lambda R_1+J_2$ is regular, then $P(\lambda)$ is regular and if $x\in\mathbb F^n\setminus\{0\}$
is an eigenvector associated with a real positive eigenvalue $\alpha$ of $P(\lambda)$ then $(\alpha J_1+J_2)x=0$.}
\item\label{sing4} {\it If the pencil $\lambda R_2+J_1$ is regular, then $P(\lambda)$ is regular and if $x\in\mathbb F^n\setminus\{0\}$
is an eigenvector associated with a real positive eigenvalue $\alpha$ of $P(\lambda)$ then $(\alpha J_1+J_2)x=0$.}
\end{enumerate}
\end{corollary}
\proof
\eqref{sing5} follows directly from~Theorem~\ref{thm:chains} and \eqref{sing0} follows from Lemma~\ref{lem:prop2} and~\eqref{sing5}.
To prove \eqref{singii} assume that $P(\lambda)$ is singular. Then for any $\lambda_0>0$ there exists a nonzero $x$ with $P(\lambda_0)x=0$
and by Lemma~\ref{lem:prop2}~i) then  $(\lambda_0 J_1+J_2)x=0$ as well. Hence $\lambda J_1+J_2$ is singular. The second claim
then follows directly. {\color{black} To see \eqref{sing3} let $P(\lambda)$ be a singular pencil, then by \eqref{sing5} the matrices $R_1$ and $J_2$ have a common kernel and the pencil $\lambda R_2+J_2$ is singular. The second statement of \eqref{sing3} follows now directly from Lemma~\ref{lem:prop2}.  The proof for~\eqref{sing4} is analogous to that for the case~\eqref{sing3}}.
\eproof
%

%

\begin{remark}\rm
The canonical forms for real or complex skew-Hermitian pencils are well-known and given in \cite{Tho91}.
These canonical forms show that in the real case all eigenvalues of a skew-symmetric matrix pencil $\lambda J_1+J_2$ necessarily
have even algebraic multiplicity. This explains why the real positive eigenvalue of the real posH pencil in Example~\ref{ex:jjb}
is a double eigenvalue.
\end{remark}

So far, we have discussed the regularity of posH pencils  as well as conditions when the spectrum
does not intersect the positive real line. Next, we will study  the index of such pencils. Although the index may be as large
as the size of the pencil, we have the following relation to the underlying skew-Hermitian matrix pencil.
\begin{theorem}\label{thm:chains2}
Let $P(\lambda)=\lambda (J_1+R_1)+(J_2+R_2)\in\mathbb F^{n,n}[\lambda]$ be a pencil of the form~\eqref{jj-pencil},
let $(x_1,\dots,x_{k})$ be a Jordan chain of length $k\geq 2$ of $P(\lambda)$ associated with the eigenvalue $\infty$ and
let $\kappa=\lfloor\frac{k+1}{2}\rfloor$.

Then we have $x_1,\dots,x_{\kappa}\in\ker R_1$ as well as $x_1,\dots,x_{\kappa-1}\in\ker R_2$
and if $k=2\kappa$ is even, then also $x_\kappa\in\ker R_2$. 
\end{theorem}
\proof
The Jordan chain $(x_1,\dots,x_{k})$ of $P(\lambda)$ associated with the eigenvalue $\infty$ satisfies
\begin{equation}\label{chain3}
(J_1+R_1)x_1=0,\quad (J_1+R_1)x_{i+1}=(J_2+R_2)x_i,\;i=1,\dots,k-1,\quad (J_2+R_2)x_{k}\neq 0
\end{equation}
or equivalently, using that $J_i^*=-J_i$ and $R_i^*=R_i$ and multiplying by $-1$,
\begin{equation}\label{chain4}
x_1^*(J_1-R_1)=0,\quad x_{i+1}^*(J_1-R_1)=x_i^*(J_2-R_2),\;i=1,\dots,k-1,\quad x_{k}^*(J_2-R_2)\neq 0.
\end{equation}
%
For the remainder of the proof we use a strategy similar to the one in the proof of Theorem~\ref{thm:chains}. From $x_1^*(J_1+R_1)x_1=0$
we get $R_1x_1=0$ and thus $J_1x_1=(J_1+R_1)x_1=0$. Furthermore, we have
\[
x_1^*(J_2+R_2)x_1=x_1^*(J_1+R_1)x_{2}=0
\]
which implies that $R_2x_1=0$.

Suppose that for some $\ell\geq 1$ with $\ell\leq\kappa-1$ we have shown $R_1x_j=0=R_2x_j$ for all $j=1,\dots,\ell$.
Then we have $\ell+2\leq\ell+2+\ell-1=2\ell+1\leq2\kappa-2+1\leq k$. Using~\eqref{chain3} and~\eqref{chain4} we obtain
\[
x_{\ell+1}^*(J_1+R_1)x_{\ell+1}=x_{\ell+1}^*(J_2+R_2)x_\ell=x_{\ell+1}^*(J_2-R_2)x_\ell=x_{\ell+2}^*(J_1-R_1)x_\ell
=x_{\ell+2}^*(J_1+R_1)x_\ell,
\]
where we have used that $R_2x_\ell=0$ and $R_1x_\ell=0$. {We repeat this procedure  $\ell-1$ times, obtaining
\[
x_{\ell+1}^*(J_1+R_1)x_{\ell+1}=x_{2\ell+1}^*(J_1+R_1)x_{1}=0
\]
which implies $R_1x_{\ell+1}=0$. If $\ell=\kappa-1$ and $k$ is odd then we are done. Otherwise (i.e. $\ell<\kappa-1$
or $k=2\kappa$ is even) we have $2\ell+2\leq k$.
Then using~\eqref{chain3},~\eqref{chain4}, $R_2x_\ell=0$ and that we just proved that
$R_1x_{\ell+1}=0$ holds, we obtain that
\begin{eqnarray*}
x_{\ell+1}^*(J_2+R_2)x_{\ell+1}&=&x_{\ell+1}^*(J_1+R_1)x_{\ell+2}=x_{\ell+1}^*(J_1-R_1)x_{\ell+2}\\
&=&x_{\ell}^*(J_2-R_2)x_{\ell+2}=x_{\ell}^*(J_2+R_2)x_{\ell+2}.
\end{eqnarray*}
As before, {\color{black} we repeat this step $\ell-1$ times, obtaining
\[
x_{\ell+1}^*(J_2+R_2)x_{\ell+1}=x_{1}^*(J_2+R_2)x_{2\ell+1}=x_{1}^*(J_1+R_1)x_{2\ell+2}=0.
\]
This implies implies $R_2x_{\ell+1}=0$.
Finally, the claim follows using an induction argument.
\eproof
\begin{corollary}\label{cor:chains2}
Let $P(\lambda)=\lambda (J_1+R_1)+(J_2+R_2)\in\mathbb F^{n,n}[\lambda]$ be a pencil of the form~\eqref{jj-pencil} and
assume that the pencil $\lambda J_1+J_2$ has at most index $\kappa$ and right minimal indices that are at most $\kappa-1$.
Then the index of $P(\lambda)$ is at most $2\kappa$.
\end{corollary}

Comparing the proofs of Theorem~\ref{thm:chains} and Theorem~\ref{thm:chains2}, we see that the main difference
is that in the proof of Theorem~\ref{thm:chains2} we can no longer use the identity $(J_2+R_2)x_k=0$, but only
$(J_1+R_1)x_1=0$. This requires us to ``push through'' the chains to the first vector $x_1$ instead of possibly to the last vector $x_k$.
This leads to the fact that not necessarily all vectors $x_1,\dots,x_k$ of the chain are in the joint kernel of the matrices
$R_1$ and $R_2$. The following examples show that the bound $\kappa$ given in Theorem~\ref{thm:chains2} is sharp.
\begin{example}\rm
Consider the pencil $P(\lambda)=\lambda (J_1+R_1)+(J_2+R_2)\in\mathbb R^{3,3}[\lambda]$ with $R_1=0$,
\[
J_1=\mat{ccc}0&0&0\\ 0&0&1\\ 0&-1&0\rix,\quad J_2=\mat{ccc}0&0&1\\ 0&0&0\\ -1&0&0\rix,\quad R_2=\mat{ccc}0&0&0\\ 0&1&0\\ 0&0&0\rix.
\]
Then $(e_1,e_2,e_3)$ is a Jordan chain of length $k=3$ of $P(\lambda)$ associated with the eigenvalue $\infty$. Here, we have
$\kappa=2$. As predicted by Theorem~\ref{thm:chains2}, we have $R_1e_1=R_2e_1=0$ and $R_1e_2=0$, but $R_2e_2\neq 0$.
The pencil $\lambda J_1+J_2$ is singular and $(e_1,e_2)$ is a singular chain of $\lambda J_1+J_2$ associated with a right
minimal index $\eta=1=\kappa-1$.
\end{example}
\begin{example}\rm
Consider the pencil $P(\lambda)=\lambda (J_1+R_1)+(J_2+R_2)\in\mathbb R^{4,4}[\lambda]$ with $R_2=0$,
\[
J_1=\mat{cccc}0&0&0&0\\ 0&0&0&1\\ 0&0&0&0\\ 0&-1&0&0\rix,\quad R_1=\mat{cccc}0&0&0&0\\ 0&0&0&0\\ 0&0&1&0\\ 0&0&0&0\rix,
\quad J_2=\mat{cccc}0&0&0&1\\ 0&0&1&0\\ 0&-1&0&0\\ -1&0&0&0\rix.
\]
Then $(e_1,e_2,-e_3,-e_4)$ is a Jordan chain of length $k=4$ of $P(\lambda)$ associated with the eigenvalue $\infty$. Again, we have
$\kappa=2$. As predicted by Theorem~\ref{thm:chains2}, we have $R_1e_1=R_2e_1=0$ and $R_1e_2=R_2e_2=0$, but $R_1e_3\neq 0$.
The pencil $\lambda J_1+J_2$ is regular and $(e_1,e_2)$ is Jordan chain of $\lambda J_1+J_2$ associated with the eigenvalue $\infty$.
\end{example}

\section{Eigenvalue localization of posH matrix pencils}\label{sec4}

In the last section we have seen that spectral properties of the underlying skew-Hermitian pencil have an important influence
on the spectral properties of a posH pencil. In view of Corollary~\ref{cor2}, one may come to the conjecture that
a posH pencil can only have eigenvalues in the right half complex plane if the underlying skew-Hermitian pencil has eigenvalues in the
right half plane or is singular. The following example shows that this conjecture is false.

\begin{example}\label{conjecture}\rm
Let
$$
J=\mat{cc} -0.1 & 1\\  0 & -0.1\rix,\quad J_1=\mat{cc} & I_2\\ -I_2  \rix, \quad
 J_2=\mat{cc} & -J\\ J^* \rix
 $$
and
$$
R_1=\mat{cccc} 1  & 1& 1 &1 \\ 1  & 1& 1 &1 \\ 1  & 1& 1 &1 \\ 1  & 1& 1 &1 \rix
\quad R_2=\mat{cccc}5 & 1& 1 &1 \\ 1  & 5& 1 &1 \\ 1  & 1& 5 &1 \\ 1  & 1& 1 &5 \rix.
$$
The eigenvalues of $P_t(\lambda)=\lambda(J_1+tR_1) + (J_2+tR_2)$ for $t\in[0,3]$ are plotted in Figure~\ref{mcdonalds}.
\begin{figure}[htb]
\includegraphics[width=400pt]{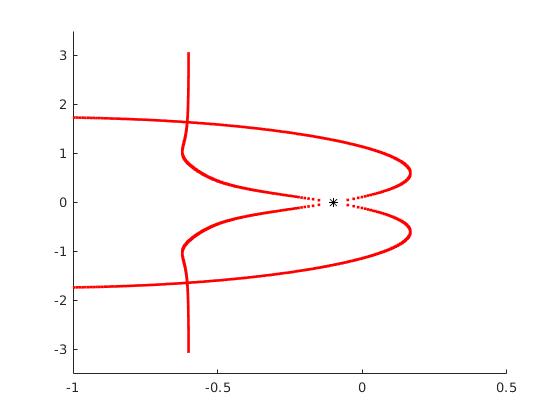}
 \caption{Figure for Example~\ref{conjecture}}\label{mcdonalds}
\end{figure}
Clearly, $P_0(\lambda)$ is a real skew-symmetric pencil, with an eigenvalue $-0.1$ (marked with star in the plot) with two
corresponding nontrivial blocks of size $2$. Even though $R_2>0$,
we observe that for some values of $t>0$ the pencil $P_t(\lambda)$ has eigenvalues in the right half plane.
\end{example}

Example~\ref{conjecture} shows that even if the pencil $\lambda J_1+J_2$ has eigenvalues in the open left half-plane, adding positive
semidefinite coefficients $R_i$ may move the eigenvalues to the right half plane. In view of this one needs further tools to localize
the eigenvalues of posH pencils. One such tool is the numerical range studied in the next subsection.

\subsection{Connections between the numerical range, common isotropic vectors and regularity of the pencil}\label{sec:commoniso}
In this subsection we employ the
 \emph{numerical range} introduced in \cite{LiR94} to obtain eigenvalue localization results.
\begin{definition}\label{def:numrange}
\rm Let $L(\lambda)\in\mathbb F^{n,n}[\lambda]$ be an $n\times n$ matrix polynomial. Then the set
\[
\mathcal W\big(L(\lambda)\big):=\big\{ \mu\in\Comp \,\big|\,   x^*L(\mu) x= 0\text{ for some }x\in\Comp^n\setminus\{ 0\}   \big\}
\]
is called the \emph{numerical range} of $L(\lambda)$.
\end{definition}

Note that in Definition~\ref{def:numrange} we take complex vectors $x\in\Comp^n$ also in the case when all coefficient are real.
We do this in order to have all finite eigenvalues of $L(\lambda)$ contained in $W(L(\lambda))$ as it is a well known fact that the
spectrum of $L(\lambda)$ is always contained in its numerical range.
This and other basic properties of the numerical range for matrix pencils and matrix polynomials are discussed in
\cite{LiR94,Psa00}. Unfortunately, there are many instances when the numerical range is the
full complex plane. First of all, this happens if $L(\lambda)$ is singular. However, this is not the only case, also for regular
pencils the numerical range can be the whole complex plane if there exist \emph{common isotropic vectors:}
\begin{definition}\rm
A matrix polynomial $L(\lambda)=\sum_{i=0}^dA_i\lambda^i\in\mathbb F^{n,n}[\lambda]$
\emph{is said to have a common isotropic vector}, 
if there exists a nonzero $x\in\mathbb C^n$ such that $x^*A_ix=0$ for all $i=0,\dots,d$.
\end{definition}
%

The following theorem relates the notion of common isotropic vectors for posH pencils to several other conditions on the pencil or
its numerical range.


%
\begin{theorem}\label{rem:nocommon}
Let $P(\lambda)=\lambda (J_1+R_1)+(J_2+R_2)\in\mathbb F^{n,n}[\lambda]$ be a posH pencil of the form~\eqref{jj-pencil},
and consider the following conditions:
\begin{enumerate}[\rm (a)]
\item $\ker R_1\cap\ker R_2=\{0\}$, \label{kerR}
\item $W\big(P(\lambda)\big) \cap (0,+\infty) = \emptyset$, \label{non+}
\item $W\big(P(\lambda)\big)\neq \Comp$, \label{nonC}
\item $P(\lambda)$ {\it has no common isotropic vector,} \label{nci}
\item $P(\lambda)$ {\it is regular.}\label{rrr}
\end{enumerate}
Then  the following implications hold:
\begin{eqnarray*}
\eqref{kerR}\To\eqref{non+}\To\eqref{nonC}\To\eqref{nci}\To\eqref{rrr},\quad \text{if } \mathbb F=\Comp\\
\eqref{kerR}\Leftrightarrow\eqref{non+}\Leftrightarrow\eqref{nonC}\Leftrightarrow\eqref{nci}\To\eqref{rrr},\quad \text{if } \mathbb F=\Real
\end{eqnarray*}
\end{theorem}
\proof
First assume that $\mathbb F=\Comp$.

\eqref{kerR}$\To$\eqref{non+}: Assume that $W\big(P(\lambda)\big)\cap\,(0,\infty)\,\neq\emptyset$. Then there exists $\alpha>0$
and a nonzero $x$ such
that $x^*W(\alpha)x=0$. Considering the real part of this equation gives $\alpha x^*R_1x+x^*R_2x=0$ which implies that $R_1x=R_2x=0$,
due to the positive semidefiniteness of $R_1$ and $R_2$. This contradicts~\eqref{kerR}.

\eqref{non+}$\To$\eqref{nonC}:
This implication is trivial.

\eqref{nonC}$\To$\eqref{nci}:
This implication is  obvious and true for arbitrary matrix polynomials, as observed in \cite{LiR94}.

\eqref{nci}$\To$\eqref{rrr}: Assume that $P(\lambda)$ is singular and take any singular chain $x_1,\dots,x_{k+1}$ as in
Theorem \ref{thm:chains}. Then, by Theorem \ref{thm:chains} we have $x_1\in\ker R_1\cap \ker R_2\cap \ker J_1$.
By Lemma \ref{lem:prop2}\eqref{c-iii} we have $x^* J_2x=0$ and consequently $x$ is a common isotropic vector.

For the case $\mathbb F=\Real$ it remains to show \eqref{nci}$\To$\eqref{kerR}:
Take a real vector $x \in\ker R_1\cap\ker R_2\setminus\{0\}$. Then clearly $x^* R_i x=x^\top J_i x=0$ for $i=1,2$, since
$x^\top J_ix$ is purely imaginary and also real. Hence, $x$ is a common isotropic vector.
\eproof

\begin{remark}\label{rem:notfull}{\rm
The following observations show the  implications between items of Theorem~\ref{rem:nocommon} that do not hold in general.

In general \eqref{rrr}$\not\To$\eqref{nci}  neither for $\mathbb F=\Comp$ nor $\mathbb F=\Real$. Indeed, take any
regular skew symmetric real pencil $\lambda J_1+J_2$ and $R_1=R_2=0$.
Then any real nonzero vector $x$ is a common isotropic vector for $P(\lambda)$.

In general \eqref{nonC}$\not\To$\eqref{non+} for $\mathbb F=\Comp$. As a counterexample take 
\begin{equation}\label{ecounter}
R_1=R_2=0,\quad J_1=\mat{cc}\ii&0\\ 0&-\ii\rix, \quad\mbox{and}\quad J_2=\mat{cc}2\ii&0\\ 0&-\ii\rix.
\end{equation}
Then  for $x=\mat{cc} x_1&x_2\rix^\top\in\mathbb C^2$ we obtain that
\[
x^*P(\alpha+\ii\beta)x=-\beta(|x_1|^2-|x_2|^2)+\ii(\alpha+2)|x_1|^2-\ii(\alpha+1)|x_2|^2
\]
which implies that $W\big(P(\lambda)\big)=\,(-\infty,-2]\cup [-1,\infty)\,$.

In general \eqref{non+}$\not\To$\eqref{kerR} for $\mathbb F=\Comp$. As a counterexample take
$J_1,J_2$ such that $\ii J_i>0$ for $i=1,2$ and $R_1=R_2=0$.
}
\end{remark}

\begin{remark}{\rm
Note that the pencil $\lambda \ii J_1+\ii J_2$ in (\ref{ecounter}) provides a
counterexample to the claim in \cite[Theorem 4.1(d)]{LiR94} that the numerical range of a Hermitian pencil equals $\mathbb R$
if both matrices are indefinite, but do not have a common isotropic vector. In fact, this was already noted in \cite{BebPNP17}.
}
\end{remark}

It remains an  open problem whether the implication \eqref{nci}$\To$\eqref{nonC} holds in the case $\mathbb F=\Comp$, but we have
the following partial result.
\begin{theorem}\label{3/4WCi} Let $P(\lambda)=\lambda (J_1+R_1)+(J_2+R_2)\in\mathbb C^{n,n}[\lambda]$ be a pencil of
the form~\eqref{jj-pencil} and let $n\geq 3$.
 If some three of the four matrices $J_1,J_2,R_1,R_2$ do not have a common isotropic vector then $W\big(P(\lambda)\big)\neq \Comp$.
\end{theorem}
\proof
Recall that, as in Lemma \ref{lem:prop2},  for $\alpha,\beta\in\Real$, it follows that
$\lambda_0=\alpha +\ii \beta\in W\big(P(\lambda)\big)$ if and only if
\begin{equation}\label{commonab}
x^* (\alpha R_1+\beta  (\ii J_1)+ R_2)x=0=x^*(\alpha (\ii J_1)- \beta  R_1+(\ii J_2))x,
\end{equation}
for some $x\neq 0$.

Assume first that the three matrices $R_1,J_1,R_2$ do not have a common isotropic vector.
Then $(0,0,0)$ does not belong to the joint numerical range
\[
W(R_1,\ii J_1,R_2)=\big\{(x^*R_1x,x^*(\ii J_1)x,x^*R_2x)\,\big|\,\norm x=1\big\}.
\]
Since $n\geq 3$, we have that $W(R_1,\ii J_1,R_2)$ is convex (see \cite{GutJK04}) and hence
$W(R_1,\ii J_1,R_2)$ coincides with its convex hull. Then \cite[Corollary 2]{Psa03} implies that
there exists nonzero scalars $\alpha,\beta,\gamma\in\Real$ such that the
the matrix $\alpha R_1+\beta \ii J_1+\gamma R_2$ is positive definite. Without loss of generality we may assume that
$\gamma=1$ (otherwise we divide by $\gamma$). Then
$\alpha R_1+\beta \ii J_1+ R_2$ is definite (positive or negative, depending on the sign of $\gamma$), and in
particular~\eqref{commonab} does not hold.
Consequently $\alpha+\ii\beta\notin W\big(P(\lambda)\big)$.

If $R_1,J_1,J_2$ do not have a common isotropic vector, then we proceed analogously, using the second equality of \eqref{commonab}.
The two other cases follow by analyzing the reversal of the pencil.
\eproof

\begin{corollary}\label{1is0}
Let $P(\lambda)$ be of the form \eqref{jj-pencil} with (at least) one of the matrices $J_1,J_2,R_1,R_2$ equal to zero, and let $n\geq 3$.
Then the implication \eqref{nci}$\To$\eqref{nonC} in Theorem~{\rm\ref{rem:nocommon}} holds.
\end{corollary}
\proof
It is enough to observe that condition \eqref{nci} is equivalent to saying that the four matrices
$J_1,J_2,R_1,R_2$ do not have a common isotropic vector. Further,  as one of them is by assumption zero,  we can apply Theorem~\ref{3/4WCi}.
\eproof

\begin{remark}{\rm
Note that statement of Theorem~\ref{3/4WCi} as well as Corollary~\ref{1is0} hold for arbitrary Hermitian matrices $R_1,R_2$
as the assumption of their nonnegativity was not used in the proof. Further, note that for $n=2$ the joint numerical range
is not necessarily convex, hence for Theorem~\ref{3/4} to hold one needs a stronger assumption. If we assume that for some
three of the four matrices $J_1,J_2,R_1,R_2$ the point $(0,0,0)$ is not in {\em the convex hull} of their numerical range,
then the proof for the case $n=2$ follows the same lines, due to Corollary 2 of \cite{Psa03}.}
\end{remark}

The following example shows that the converse statement in Theorem~\ref{3/4WCi} does not hold.

\begin{example}{\rm
Let
\[
P(\lambda)=\mat{cccc} 1 \\ & \lambda \\ && \ii\\ &&& \lambda\ii\rix
=\lambda\mat{cccc}0\\ &1\\ &&0\\ &&&\ii\rix+\mat{cccc}1\\ &0\\ &&\ii\\ &&&0\rix.
\]
Then each three the four matrices $J_1,J_2,R_1,R_2$ have a common kernel and consequently a common isotropic vector. However, the
numerical range is contained in the left half plane. Indeed, from $x^* P(\alpha+i\beta)x=0$ we obtain
\[
\alpha|x_2|^2+|x_1|^2-\beta |x_4|^2=0\quad\mbox{and}\quad \beta |x_2|^2+\alpha |x_4|^2+|x_3|^2=0.
\]
If $\alpha>0$, then the first equality implies $x_1=0$ and $\beta\geq 0$, but then the second equality gives
$x_2=x_3=x_4=0$. Hence, if $\alpha>0$, then $\alpha+i\beta$ is not in the numerical range of $P(\lambda)$ for any $\beta$ .
}
\end{example}

\subsection{Localizing the numerical range in a pacman-like shape}\label{sec:pack}
In this subsection we present a localization result for the numerical range of posH pencils.
For this we introduce the notation
\begin{equation}\label{betapmdef}
\beta_{\pm}:=\sup\set{ \beta \geq 0\,|\,R_1 + R_2 \pm \beta (\ii J_1) >0 },
\end{equation}
provided that the corresponding set under the supremum is nonempty (otherwise we do not define the symbol).
Both $\beta_{\pm}$ are well-defined and are either positive or equal to $\infty$ if (but not only if) $\ker R_1\cap \ker R_2 =\set 0$.
In such case
\begin{equation}\label{ebound}
\beta_{\pm}  \geq \sigma_{\min}(R_1+R_2)/\norm{J_1}.
\end{equation}
These bounds may be, however, far from optimal and a direct numerical estimation may give better bounds.


\begin{theorem}\label{3/4} Let $P(\lambda)=\lambda (J_1+R_1)+(J_2+R_2)\in\mathbb C^{n,n}[\lambda]$ be a pencil of the form~\eqref{jj-pencil}
and let $n\geq 3$.
\begin{enumerate}[\rm (i)]
\item\label{WCii} If $R_1,R_2,J_1$ do not have a common isotropic vector then either $\beta_+>0$ and
\begin{equation}\label{beta+}
W\big(P(\lambda)\big)\,\cap\,\big\{z\in\Comp\,\big|\, \RE z>0,\ 0 \leq \IM z < \beta_+ ,\ 0\leq \arg z< \arctan(\beta_+)\big\}=\emptyset,
\end{equation}
or $\beta_->0$ and
\begin{equation}\label{beta-}
W\big(P(\lambda)\big)\,\cap\,\big\{z\in\Comp\,\big|\, \RE z>0,\ - \beta_- < \IM z \leq 0 ,\ -\arctan(\beta_-) < \arg z\leq  0\big\}=\emptyset,
\end{equation}
where we use the convention $\arctan(\infty)=\frac{\pi}{2}$.
\item \label{WCiii} If  $R_1,R_2$ do not have a common isotropic vector, i.e.~if $\,\ker R_1\cap\,\ker R_2=\{0\}$ then
 $\beta_-,\beta_+>0$ and both \eqref{beta+} and \eqref{beta-} hold.

\end{enumerate}
\end{theorem}
\proof
\eqref{WCii} Assume that $R_1,R_2,J_1$ do not have a common isotropic vector. As in the proof of Theorem~\ref{3/4WCi} we get
that there exist nonzero $\alpha,\beta,\gamma$ such that  the matrix $\alpha R_1+\beta \ii J_1+\gamma R_2$ is positive definite. {
Due to the positive semidefiniteness of $R_1$ and $R_2$, we can increase $\alpha$ and $\gamma$ without changing the property of
positive definiteness, and hence we can assume $\alpha=\gamma>0$. Dividing by $\alpha$, we obtain $ R_1+\beta_0 \ii J_1+ R_2>0$
for some $\beta_0\in\Real$. Assume first that $\beta_0>0$. By the positive semidefiniteness of $R_1$ and $R_2$, we have that
\[
a R_1 + b (\ii J_1) + c R_2 >0,\quad a,c\geq 1,\ 0\leq b\leq\beta_0,
\]
consequently $\beta_+\geq\beta_0>0$. Hence,
\[
x^*\left( \frac ac R_1 + \frac bc (\ii J_1) +  R_2\right)x  >0,\quad a,c\geq 1,\ 0\leq b<\beta_+,\ x\neq 0.
\]
This means that for all  $a,c\geq 1,\ 0\leq b<\beta_+$, the point $\lambda_0=\frac ac+\ii \frac bc$ is not in $W(P(\lambda))$,
see once again \eqref{commonab}.
It is then an easy calculation to see that \eqref{beta+} holds: take $\lambda_0$ in the set excluded from the numerical range
by~\eqref{beta+}, i.e., $\RE\lambda_0>0$, $0\leq \IM \lambda_0<\beta_+$ and $\IM \lambda_0<\beta_+\RE\lambda_0$. If $\RE\lambda_0\geq 1$
then we take $0\leq b<\beta_+$ and $c\geq 1$ such that $\IM\lambda_0=\frac bc$ and set $a=c\ \RE\lambda_0$, which is necessarily
greater or equal to one. If  $\RE\lambda_0< 1$ then we take $a=1$, $c=\frac 1{{\rm Re}\lambda_0}>1$, and since
$\IM\lambda_0<\beta_+\RE\lambda_0=\frac{\beta_+}c$ we can find $0\leq b<\beta_+$ such that $\IM\lambda_0=\frac bc$.

Analogously, if $\beta_0<0$ then $\beta_- \geq -\beta_0>0$ and we obtain \eqref{beta-}.

Finally, if  $\beta_0=0$ then we are in the case $R_1+R_2>0$, described by statement \eqref{WCiii}. In such a situation we may find
(by continuity) $\beta_0'>0$ as well as $\beta_0'<0$ such that $ R_1+\beta_0' \ii J_1+ R_2>0$. Consequently, $\beta_\pm>0$ and the
proof follows the same lines as before.  }
%
\eproof
Let us illustrate the Theorem~\ref{3/4} with an example.
\begin{example}\label{ex1}{\rm
Consider a $10\times10$ pencil, generated randomly by \texttt{matlab} via
\begin{verbatim}
R1=randn(10); R1=0.04*R1'*R1;
R2=randn(10); R2=0.04*R2'*R2;
J1=rand(10); J1=J1-J1';
J2=rand(10); J2=J2-J2';
\end{verbatim}
In our particular example $R_1$ and $R_2$ are close to singular (their smallest eigenvalues are of order $10^{-3}$ and
$10^{-5}$ respectively), but thee sum has the smallest eigenvalue of order $10^{-1}$.
The approximation of the numerical range is plotted in blue using $10^5$ random points, the eigenvalues are marked with red circles,
see Figure~\ref{f2}.
We numerically calculate $\beta_{\pm}=0.2032$, and since  all matrices are real, the two values  coincide.
For comparison, $\sigma_{\min}(R_1+R_2)/\norm{J_1}=0.0588$, cf. \eqref{ebound}. The set excluded from the numerical
range due to Theorem~\ref{3/4}\eqref{WCiii} is displayed between the blue and orange  line.
\begin{figure}[htb]
\begin{center}
\includegraphics[width=300pt]{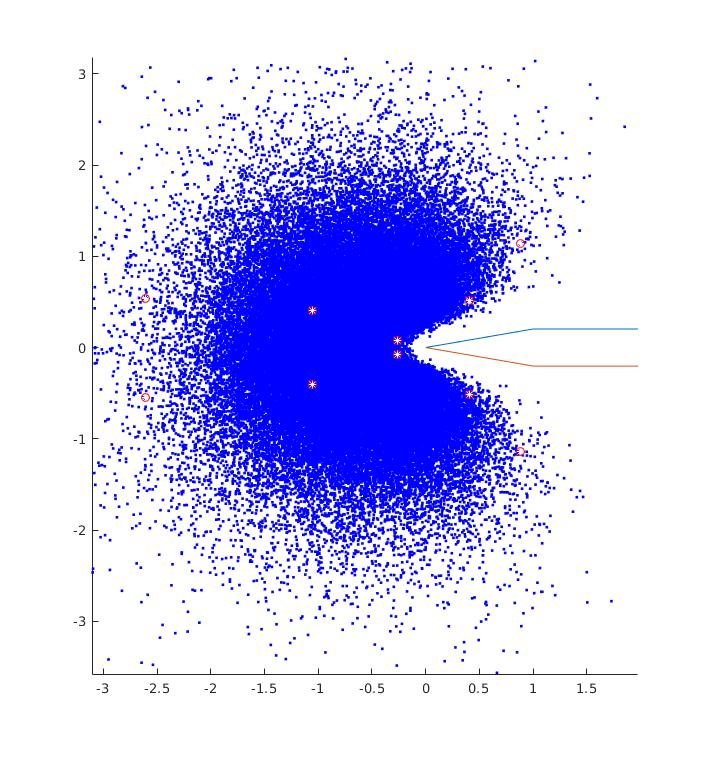}
\end{center}
\caption{Numerical range and spectrum of the pencil in Example~\ref{ex1}.}\label{f2}
\end{figure}
}
\end{example}
\begin{remark}{\rm
The border of the excluded region in \eqref{beta+} is a line, that splits up at a point $1+\ii \beta_+$. One can easily
generate a different splitting point by introducing
\[
\beta_{\pm}(t):=\sup\set{ \beta \geq 0\,|\,tR_1 + R_2 \pm \beta (\ii J_1) >0 },\quad t>0,
\]
obtaining
\[
W\big(P(\lambda)\big)\,\cap\,\big\{z\in\Comp\,\big|\, \RE z>0,\ 0 \leq \IM z < \beta_+(t) ,\ 0\leq \arg z<
\arctan\big(\beta_+(t)/t\big)\big\}=\emptyset.
\]
Further, if $R_2$ is invertible and $J_1\neq 0$, then $\beta_+(t)\geq\sigma_{\min}(R_2)/\norm{J_1} $.
Letting $t\to 0$ we obtain  that a strip is excluded from the numerical range:
\[
W\big(P(\lambda)\big)\,\cap\,\big\{z\in\Comp\,\big|\, \RE z>0,\ 0 \leq \IM z < \sigma_{\min}(R_2)/\norm{J_1}  \big\}=\emptyset,
\]
which is also visible in Figure~\ref{f2}.}
\end{remark}
\subsection{Localizing the numerical range and the spectrum in the left half-plane}\label{sec:left}

In this subsection we investigate sufficient conditions for the numerical range and the spectrum of pencils $P(\lambda)$
of the form \eqref{jj-pencil} to be contained in the left half plane.
\begin{lemma}\label{lem:new-remodelled-twice}
Let $P(\lambda)=\lambda (J_1+R_1)+(J_2+R_2)\in\mathbb C^{n,n}[\lambda]$ be a pencil of the form~\eqref{jj-pencil}.
%
If $\lambda_0\in\mathbb C$ and $x\in\mathbb C^n$ are such that $x^*P(\lambda_0)x=0$, $x^*(J_1+R_1)x\neq 0$ and
\begin{equation}\label{EE-JJ}
-x^*R_1xx^*R_2x+x^*J_1xx^*J_2x\leq 0,
\end{equation}
then $\RE\lambda_0\leq0$.
\end{lemma}
\proof
The proof follows from the fact that $\lambda_0 x^*(J_1 +R_1)x =- x^*(J_2+R_2)x$ implies that
\[
\RE\lambda_0=\frac{\RE\big(\!- x^*(J_2+R_2)x\cdot\overline{x^*(J_1 +R_1)x}\,\big)}{|   x^*(J_1+R_1)x |^2}=
 \frac{-x^*R_1xx^*R_2x+x^*J_1xx^*J_2x}{|   x^*(J_1+R_1)x |^2}\leq 0.\quad\Box
\]

\begin{theorem}\label{thm:lefthalfplane}
Let $P(\lambda)=\lambda (J_1+R_1)+(J_2+R_2)\in\mathbb C^{n,n}[\lambda]$ be a posH pencil, i.e., of the form \eqref{jj-pencil}, such that
\begin{equation}\label{EE-JJx}
-x^*R_1xx^*R_2x+x^*J_1xx^*J_2x\leq 0,\quad x\in\Comp^n
\end{equation}
and (at least) one of the following two conditions,  cf.~{\rm Theorem~\ref{rem:nocommon}}, hold.
%
\begin{enumerate}[\rm (a)]
\item\label{t-i}  {\it $P(\lambda)$ has no common isotropic vector.}
\item\label{t-ii} {\it $P(\lambda)$ is a regular pencil, and $W(\lambda J_1+J_2)$ is contained in the closed left half-plane.}
\end{enumerate}
%
Then the numerical range of $P(\lambda)$, and thus all finite eigenvalues, are contained in the closed left half-plane.
\end{theorem}
\proof
First assume  that \eqref{t-i} holds and take $\lambda_0\in \mathcal W\big(P(\lambda)\big)$, i.e., for some $x\neq 0$ we have
$
\lambda_0 x^*(J_1 +R_1)x =- x^*(J_2+R_2)x.
$
Note that due to \eqref{t-i} we necessarily have  $x^*(J_1+R_1)x\neq 0$ and Lemma \ref{lem:new-remodelled-twice} can be applied.

Secondly assume \eqref{t-ii} and let $x\neq 0$ be such that  $x^*P(\lambda_0)x=0$. If $x^*(J_1+R_1)x\neq 0$, then the claim follows
by Lemma~\ref{lem:new-remodelled-twice}. If $x^*(J_1+R_1)x= 0$,  then by Lemma~\ref{lem:prop2} we have $R_1x=R_2x=0$ and
$x^*(\lambda_0 J_1+J_2)x=0$. Since the numerical range of the pencil $\lambda J_1+J_2$ is contained in the closed left half plane, we
obtain that $\RE\lambda_0\leq 0$.
\eproof
\begin{corollary} Let $P(\lambda)=\lambda (J_1+R_1)+(J_2+R_2)\in\mathbb C^{n,n}[\lambda]$ be a posH pencil, i.e., of the form~\eqref{jj-pencil},
such that $\ker R_1\cap \ker R_2=\set 0$ and \eqref{EE-JJx} holds. Then $P(\lambda)$ is regular with the numerical range
(and hence all finite eigenvalues)  contained in the closed left half-plane.
\end{corollary}

Condition \eqref{t-ii} in Theorem~\ref{thm:lefthalfplane} is rather strong. Relaxing it, we are still able to make a statement
on the spectrum. 
\begin{theorem}\label{thm:lefthalfplane2}
Let $P(\lambda)=\lambda (J_1+R_1)+(J_2+R_2)\in\mathbb C^{n,n}[\lambda]$ be a posH pencil, i.e., of the form~\eqref{jj-pencil},
such that \eqref{EE-JJx} holds.
If $P(\lambda)$ is a regular pencil,  $\lambda J_1+J_2$ has all minimal indices equal to zero (if there are any), and all its
eigenvalues in the closed left half-plane, then  all the eigenvalues of $P(\lambda)$ are contained in the closed left half-plane.
\end{theorem}
\proof
Let $\lambda_0$ be an eigenvalue of $P(\lambda)$ with corresponding eigenvector $x$.
If $x^*(J_1+R_1)x\neq 0$ then the result follows from Lemma~\ref{lem:new-remodelled-twice}. If $x^*(J_1+R_1)x= 0$,  then by
Lemma~\ref{lem:prop2} we have $R_1x=R_2x=0$ and hence $P(\lambda_0)x=0$ implies $(\lambda_0 J_1+J_2)x=0$. If $\RE\lambda_0>0$
then by assumption $\lambda_0$ is not an eigenvalue of $\lambda_0 J_1+J_2$, and hence, given that all minimal indices
of $\lambda J_1+J_2$ are
equal to zero, $x$ is in the joint kernels of $J_1$, $J_2$, $R_1$, and $R_2$. Therefore $P(\lambda)$ is singular, which is a
contradiction.
\eproof
%

%
%
%
We see that condition~\eqref{EE-JJx} plays a crucial role in the characterization of stable pencils of the form~\eqref{jj-pencil}.
Unfortunately, this condition is in general hard to check. For this reason we present a result under stronger assumptions that can
in general be verified more easily. Let $A\otimes B$ denote the Kronecker product of matrices $A$ and $B$,
see e.g.~\cite{LieM15a}. Note that $J_1\otimes J_2$ is a Hermitian matrix if $J_1$ and $J_2$ are skew-Hermitian.

\begin{proposition}\label{prop:kronecker}
Let $P(\lambda)=\lambda (J_1+R_1)+(J_2+R_2)\in\mathbb C^{n,n}[\lambda]$ be a  posH pencil, i.e., of the form~\eqref{jj-pencil}.
\begin{enumerate}[\rm (i)]
\item \label{EE0} If $\lambda_{\min}(R_1)\lambda_{\min}(R_2)\geq\norm{J_1}\norm{J_2}$, where $\lambda_{\min}$ denotes the smallest
eigenvalue, then \eqref{EE-JJx} holds.
\item\label{EEi} If $J_1\otimes J_2-R_1\otimes R_2 \leq 0$, then \eqref{EE-JJx} holds.
\item\label{EEii} If all finite eigenvalues of $\lambda J_1+J_2$ are real, nonpositive and semisimple,  if the eigenvalue $\infty$
(if it exists) is semisimple, and if all  minimal indices are zero (if there are any), then
\[
x^*J_1xx^*J_2x\leq 0
\]
holds for all $x\in\mathbb C^n$ and thus, in particular, \eqref{EE-JJx} holds.
\item\label{EE3} If $J_1,J_2,R_1,R_2\in\Real^n$ 
 then \eqref{EE-JJx} is equivalent to
\begin{equation}\label{EE-JJReal}
4(\xi^\top J_1\eta)(\xi^\top J_2\eta)\leq(\xi^\top R_1\xi+\eta^\top R_1\eta)(\xi^\top R_2\xi+\eta^\top R_2\eta),\quad \xi,\eta\in\Real^n.
\end{equation}
\end{enumerate}
\end{proposition}
\proof
\eqref{EE0} This follows, since  $x^*R_i x \geq \lambda_{\min}$ for every vector $x$.\\
\eqref{EEi} Since by well-known properties on the Kronecker product, see e.g.~\cite{LieM15a}, condition~\eqref{EE-JJx} can be rewritten as
\[
(x\otimes x)^*(J_1\otimes J_2)(x\otimes x)-(x\otimes x)^*(R_1\otimes R_2)(x\otimes x)\leq 0,
\]
the assertion follows.

\eqref{EEii}
Considering the Hermitian pencil $\lambda i J_1 +i J_2$, we may assume that this pencil is in
the Hermitian canonical form of \cite{Tho76}. By the assumptions on the spectrum of $\lambda J_1+J_2$ (and thus
$\lambda i J_1 +i J_2$), if follows that both $i J_1=\diag(a_1,\dots,a_n)$ and $i J_2=\diag(b_1,\dots,b_n)$ are diagonal
and that two diagonal elements in the same position are either both nonnegative or both nonpositive, i.e., $a_ib_i\geq 0$
for $i=1,\dots,n$. Then for any vector $x=[x_1,\dots,x_n]^\top\in\mathbb C^n$ we have
$\overline{x}_ja_ix_j\overline{x}_jb_ix_j=a_ib_i|x_j|^4\geq 0$, which implies that $x^*i J_1xx^*i J_2x \geq 0$, or, equivalently,
$x^* J_1xx^* J_2x \leq 0$.

\eqref{EE3} It is enough to take $x=\xi+i \eta$ in \eqref{EE-JJx}  to see the equivalence.
\eproof

%
%
%

Although the characterization of the spectrum using the numerical range and the existence of common isotropic vectors is not complete
(the question whether the implication \eqref{nci}$\To$\eqref{nonC} holds in Theorem~\ref{rem:nocommon} is still open), it is
rather surprising to observe how many properties carry over from dH to posH pencils.
In the next section we discuss the extension of some of these results to matrix polynomials.

\section{Matrix polynomials with positive semidefinite Hermitian coefficients}\label{sec5}

In this section we investigate matrix polynomials with positive semidefinite Hermitian coefficients
\begin{equation}\label{PSC}
P(\lambda)=\lambda^d A_d+\cdots+\lambda A_1+A_0,\quad  d\geq 1,\quad A_j^*=A_j\geq 0,\ j=0,\dots,d.
\end{equation}
The analysis when such a polynomial is singular was presented in \cite{MehMW21}, where it was shown that all left
and right minimal indices cannot exceed zero. Regarding the spectrum,
it is well known that the eigenvalues are always in the closed left half plane if the degree of the polynomial
is less than two, see Theorem~\ref{stable=dh} and \cite[Corollary 4.9]{MehMW18}. Unfortunately, this is no longer true if the degree
exceeds three. As an example consider the scalar polynomial $P(\lambda)=\lambda^3+1$ that has eigenvalues in the open right half plane.

Due to the observation on the spectrum of the previous paragraph, it is clear that matrix polynomials as in~\eqref{PSC} can in
general not be linearized by dh pencils. Instead, we will show that they can be linearized by posH pencils.
In the following, we will first give bounds for the index of the polynomial using the results of Section~\ref{sec3}.
Next, we will localize the spectrum using Section~\ref{sec4}.

\subsection{The index of a matrix polynomial with positive semidefinite Hermitian coefficients}

Possible linearizations of matrix polynomials in \eqref{PSC} leading to posH matrix pencils have been derived in the literature.
First, let us assume that the degree $d$ of the matrix polynomial is odd, then we can find a posH linearization via a block
symmetric linearization presented in \cite{AntV04} that was later identified as a special instance of a generalized Fiedler
pencil and hence is a strong linearization of the given matrix polynomial \cite{DeTDM10}.
\begin{remark}\label{rem:5.4.21}{\rm
Let $d=2\delta-1$ be odd and let $P(\lambda)=\sum_{i=0}^dA_i\lambda^i\in\mathbb F^{n,n}[\lambda]$ with $A_i^*=A_i\geq 0$ for
$i=0,\dots,d$. Then by \cite{AntV04} the block $dn\times dn$ matrix pencil
$S_P(\lambda)=\mat{c}S_{ij}(\lambda)\rix_{i,j=1,\dots,d}$ with
\begin{eqnarray*}
S_{2j-1,2j-1}(\lambda)&=&\lambda A_{2j-1}-A_{2j-2},\quad j=1,\dots,\delta \\
S_{2j-1,2j}(\lambda)&=&S_{2j,2j-1}(\lambda)=\lambda I_n,\quad j=1,\dots,\delta-1\\
S_{2j,2j+1}(\lambda)&=&S_{2j+1,2j}(\lambda)=I_n,\quad j=1,\dots,\delta-1
\end{eqnarray*}
and $S_{ij}(\lambda)=0$ for all remaining blocks is a linearization (in fact,
by \cite{DeTDM10} a strong linearization) of the matrix polynomial $P(\lambda)$.
Multiplying the $2j$-th block row with $-1$ for $j=1,\dots,\delta-1$, we obtain a strong linearization
$S_P(\lambda)=\lambda (J_1+R_1)+(J_2+R_2)$ of $P(\lambda)$ having coefficients with positive semidefinite Hermitian parts of the form
\begin{eqnarray}
J_1&=&\diag\left(\mat{cc}0&I_n\\ -I_n&0\rix,\dots,\mat{cc}0&I_n\\ -I_n&0\rix,0\right),\label{5.4.21.1}\\
J_2&=&\diag\left(0,\mat{cc}0&-I_n\\ I_n&0\rix,\dots,\mat{cc}0&-I_n\\ I_n&0\rix\right),\label{5.4.21.2}\\
R_1&=&\diag\left(\mat{cc}A_1&0\\ 0&0\rix,\dots,\mat{cc}A_{d-2}&0\\ 0&0\rix,A_d\right),\label{5.4.21.3}\\
R_2&=&\diag\left(\mat{cc}A_0&0\\ 0&0\rix,\dots,\mat{cc}A_{d-3}&0\\ 0&0\rix,A_{d-1}\right),\label{5.4.21.4}
\end{eqnarray}
where each $0$ stands for the $n\times n$ zero matrix.
}
\end{remark}

For the special cases $d=3$ and $d=5$ we have the following pencils $S_P(\lambda)$.
\begin{example}{\rm
Let $A_i\in\mathbb F^{n,n}$ and $A_i^*=A_i\geq 0$ for $i=0,\dots,5$ and consider the matrix polynomials
\begin{eqnarray*}
P(\lambda)&=&A_3\lambda^3+A_2\lambda^2+ A_1\lambda +A_0,\\
Q(\lambda)&=&A_5\lambda^5+A_4\lambda^4+A_3\lambda^3+A_2\lambda^2+A_1\lambda +A_0.
\end{eqnarray*}
Then we have the linearizations
\[
S_P(\lambda)=\lambda\mat{ccc}A_1 & I_n & 0\\ -I_n & 0& 0\\ 0&0&A_3\rix+\mat{ccc}A_0&0&0\\ 0&0&-I_n\\ 0&I_n&A_2\rix
\]
and
\[
S_Q(\lambda)=\lambda\mat{ccccc}A_1&I_n&0&0&0\\ -I_n&0&0&0&0\\ 0&0&A_3&I_n&0\\ 0&0&-I_n&0&0\\ 0&0&0&0&A_5\rix+
\mat{ccccc}A_0&0&0&0&0\\ 0&0&-I_n&0&0\\ 0&I_n&A_2&0&0\\ 0&0&0&0&-I_n\\ 0&0&0&I_n&A_4\rix.
\]
}
\end{example}

If, on the other hand, the degree $d$ of a  matrix polynomial in \eqref{PSC} is even, then a similar block symmetric linearization
is only known
for the case that one of the coefficients $A_0$ or $A_d$ is invertible, see again \cite{AntV04}. Since we are particularly interested
in the index of the pencil, we focus on the case that $A_0$ is invertible.

\begin{remark}\label{rem:5.4.21.2}{\rm
Let $d=2\delta$ be even and let $P(\lambda)=\sum_{i=0}^dA_i\lambda^i\in\mathbb F^{n,n}[\lambda]$ with $A_i^*=A_i\geq 0$ for
$i=1,\dots,d$ and $A_0^*=A_0>0$. Then by \cite{MacMMM09b} the block $dn\times dn$ matrix pencil
$S_P(\lambda)=\mat{c}S_{ij}(\lambda)\rix_{i,j=1,\dots,d}$
with
\begin{eqnarray*}
S_{1,1}(\lambda)&=&A_{0}^{-1},\\
S_{2j,2j}(\lambda)&=&\lambda A_{2j-1}-A_{2j},\quad j=1,\dots,\delta, \\
S_{2j-1,2j}(\lambda)&=&S_{2j,2j-1}(\lambda)=\lambda I_n,\quad j=1,\dots,\delta,\\
S_{2j,2j+1}(\lambda)&=&S_{2j+1,2j}(\lambda)=I_n,\quad j=1,\dots,\delta-1,
\end{eqnarray*}
and $S_{ij}(\lambda)=0$ for all remaining blocks, is a linearization (in fact,
by \cite{DeTDM10} a strong linearization) of the matrix polynomial $P(\lambda)$.
Multiplying the $(2j-1)$-th block row with $-1$ for $j=1,\dots,\delta-1$, we obtain a strong linearization
$S_P(\lambda)=\lambda (J_1+R_1)+(J_2+R_2)$ of $P(\lambda)$ having coefficients with positive semidefinite Hermitian parts of the form
\begin{eqnarray}
J_1&=&\diag\left(0,\mat{cc}0&-I_n\\ I_n&0\rix,\dots,\mat{cc}0&-I_n\\ I_n&0\rix,0\right),\label{5.4.21.5}\\
J_2&=&\diag\left(\mat{cc}0&I_n\\ -I_n&0\rix,\dots,\mat{cc}0&I_n\\ -I_n&0\rix\right),\label{5.4.21.6}\\
R_1&=&\diag\left(A_0^{-1},\mat{cc}A_2&0\\ 0&0\rix,\dots,\mat{cc}A_{d-2}&0\\ 0&0\rix,A_{d}\right),\label{5.4.21.7}\\
R_2&=&\diag\left(\mat{cc}0&0\\ 0&A_1\rix,\dots,\mat{cc}0&0\\ 0&A_{d-1}\rix\right).\label{5.4.21.8}
\end{eqnarray}
%
}
\end{remark}
For the special cases $d=4$ we have the following pencil $S_P(\lambda)$.
\begin{example}\rm
Let $A_i\in\mathbb F^{n,n}$ and $A_i^*=A_i\geq 0$ for $i=0,\dots,4$ and consider the matrix polynomial
$P(\lambda)=A_4\lambda^4+A_3\lambda^3+A_2\lambda^2+ A_1\lambda+A_0.$ Then we have
\[
S_P(\lambda)=\lambda\mat{cccc}A_0^{-1}&0&0&0\\ 0&A_2&-I_n&0\\ 0&I_n&0&0\\ 0&0&0&A_4\rix
-\mat{cccc}0 & I_n & 0&0\\ -I_n & A_1& 0&0\\ 0&0&0&I_n\\ 0&0&-I_n&A_3\rix.
\]
\end{example}
Concerning the possible index of matrix polynomials with Hermitian positive semidefinite coefficients, we are able to make general
statements using the results from Section~\ref{sec3}.

\begin{theorem}\label{thm:indexd}
Let $P(\lambda)=\sum_{i=0}^d\lambda^iA_i\in\mathbb F^{n,n}[\lambda]$ with $A_i^*=A_i\geq 0$ for
$i=0,\dots,d$, where we assume that $A_0$ is invertible if $d$ is even.
Then the index of $P(\lambda)$ is at most $d$.
\end{theorem}
\proof
First, let $d=2\delta-1$ be odd and let $S_P(\lambda)=\lambda (J_1+R_1)+(J_2+R_2)$ be the linearization of $P(\lambda)$ from
Remark~\ref{rem:5.4.21}.
Let $k$ be the index of $S_P(\lambda)$ (and thus also of $P(\lambda)$) and let $(x_1,\dots,x_k)$ be a Jordan chain of
$S_P(\lambda)$ associated with the eigenvalue $\infty$.
Suppose that $k>d$, i.e., $k\geq 2\delta$. 
By Theorem~\ref{thm:chains2} we have that
\[
x_1,\dots,x_\delta\in\ker R_1\cap\ker R_2.
\]
In particular, this implies $J_1x_1=0$ and $J_1x_{j+1}=J_2x_j$ for
$j=1,\dots,\delta-1$. By~\eqref{5.4.21.1} it follows that
$x_1$  must be of the form $\mat{cccc}0 &\dots &0 &x^T\rix^T$ with some $x\neq 0$ and $x_j$
 must be the vector
that has $x$ in the $(2(\delta-j)+1)$-th block component and is zero anywhere else for $j=1,\dots,\delta$. Hence, $J_2x_\delta=0$,
and since also $R_2x_\delta=0$, we find that $(x_1,\dots,x_k)$ is not a Jordan chain of $S_P(\lambda)$ associated with the
eigenvalue $\infty$, which is a contradiction.
%
%
%
%

Secondly, let $d=2\delta$ be even and let $S_P(\lambda)=\lambda (J_1+R_1)+(J_2+R_2)$ be the linearization of $P(\lambda)$ from
Remark~\ref{rem:5.4.21.2}. Let $k$ be the index of $S_P(\lambda)$ (and thus also of $P(\lambda)$) and let $(x_1,\dots,x_k)$
be a Jordan chain of $S_P(\lambda)$ associated with the eigenvalue $\infty$.

Suppose that $k>d$, i.e., $k\geq 2\delta+1$. 
By Theorem~\ref{thm:chains2} we have that
$$
x_1,\dots,x_\delta\in\ker R_1\cap\ker R_2,\quad x_{\delta+1}\in\ker R_1.
$$
In particular, $J_1x_1=0$ and $J_1x_{j+1}=J_2x_j$ for
$j=1,\dots,\delta-1$. By~\eqref{5.4.21.5} and~\eqref{5.4.21.7} if follows that
$x_1$ must be of the form $\mat{cccc}0&\dots&0&x^T\rix^T$ with some $x\neq0$. Then it follows from~\eqref{5.4.21.5} and~\eqref{5.4.21.6}
that $x_j$ must be the vector
that has $x$ in the $(2(\delta-j)+2)$-th block component and is zero anywhere else for $j=1,\dots,\delta$. Furthermore,
we have that $J_2x_\delta$ is not in the range of $J_1$. However, since $(x_1,\dots,x_k)$ is a Jordan chain of $S_P(\lambda)$
associated with the eigenvalue $\infty$, we have that
\[
J_2x_\delta =(R_2+J_2)x_\delta= (J_1+R_1)x_{\delta+1}=J_1x_{\delta+1},
\]
which again is a contradiction.
%
\eproof
The bound in Theorem~\ref{thm:indexd} is sharp as the following example shows.
\begin{example}{\rm
Consider the (scalar) $1\times 1$ matrix polynomial $P(\lambda)=\sum_{i=0}^da_i\lambda^i$ with $a_0=1$ and $a_1=\cdots=a_d=0$.
It is easy to check that the chain $(e_d,e_{d-2},\dots,e_1,e_2,e_4,\dots,e_{d-1})$ is a Jordan chain of the
pencil $S_P(\lambda)$ as in Remark~\ref{rem:5.4.21} associated with the eigenvalue $\infty$ if $d$ is odd.
If, on the other hand, $d$ is even, then
the chain $(e_d,e_{d-2},\dots,e_2,e_1,e_3,\dots,e_{d-1})$ is a Jordan chain of the
pencil $S_P(\lambda)$ as in Remark~\ref{rem:5.4.21.2} associated with the eigenvalue $\infty$.
Thus, in both cases we find that $P(\lambda)$ has index $d$.
}
\end{example}

\subsection{Eigenvalue localization for matrix polynomials with positive semidefinite Hermitian coefficients}
In this section we present eigenvalue localization results for matrix polynomials with positive semidefinite Hermitian coefficients.
First let us observe that the spectrum of such matrix polynomials is still restricted, due to an analogous result on scalar
polynomials with positive coefficients.

%
%
\begin{theorem}\label{thm:polys}
Let $P(\lambda)=\lambda^d A_d+\cdots+\lambda A_1+A_0$, $d\geq 1$ be a regular $n\times n$ complex matrix
polynomial, where $A_j^*=A_j\geq 0$ for $j=0,\dots,d$.
%
Then the numerical range $\mathcal W\big(P(\lambda)\big)$ and hence the spectrum of $P(\lambda)$
is contained in $\big\{z\in\mathbb C\,:\,|\arg(z)|\geq\frac{\pi}{d}\big\}\cup\{0\}$.
\end{theorem}
\proof
 Let $\mu\in \mathcal W\big(P(\lambda)\big)$. Then there exists $x\in\mathbb C^n\setminus\{0\}$ such that $\mu$ is a root of
the polynomial $p(\lambda):=x^*P(\lambda)x=a_d\lambda^d +\cdots+a_1\lambda +a_0$ with $a_i=x^*A_ix\geq 0$ for $i=0,\dots,d$.
If $p$ is the zero polynomial, then $A_dx=\cdots=A_0x=0$ and hence $P(\lambda)$ is singular, which is a contradiction.
Hence there exists an index $i\in \{0,\ldots,d\}$ with $a_i\neq 0$. Let $\ell$ and $m$ be the maximal and minimal indices $i$
such that $a_i\neq 0$,
respectively. If $\ell=m$, then $p(\lambda)=a_m\lambda^m$ only has zero as its root which implies $\mu=0$.
Otherwise, we have $p(\lambda)=\lambda^m\widetilde p(\lambda)$, where $\widetilde p(\lambda)$ is a
polynomial of degree $\ell-m\geq 1$ with nonnegative coefficients and with the leading and last term being additionally nonzero.
If $\ell-m=1$, then $\widetilde p(\lambda)$ necessarily has a negative root and the result is trivial. If $\ell-m\geq 2$, then by
Theorem 4.1 of \cite{CowT54} the polynomial $\widetilde{p}(\lambda)$
has all its root outside the angle $|\arg(\lambda)|<\frac{\pi}{\ell-m}$ and the claim follows. \eproof



Theorem~\ref{thm:polys} shows that for matrix polynomials with Hermitian positive semidefinite coefficients and degree at least three,
the spectrum is not automatically contained in the closed left half plane and therefore, it is necessary to decide whether
this is the case or not. This can be done by applying the results on posH pencils from the previous sections on appropriate
linearizations for the given matrix polynomials. In the following, we explicitly reformulate some results in terms of the
coefficients of the matrix polynomial for the important special case $d=3$ by using a special posH linearization that only
contains the coefficients of the matrix polynomial as nonzero blocks.

%
Let $L(\lambda)=\lambda^3 A_3+\lambda^2 A_2+\lambda A_1 +A_0$ be a complex matrix polynomial with $A_0,A_3>0$.
Due to \cite{HigMMT06b} we have the following strong linearization
\[
\lambda \mat{ccc}0&\!\!A_3\!\!&0\\ \!A_3\!&A_2&0\\0&0&-A_0\rix\!+\mat{ccc}-A_3&0&0\\ 0&A_1&\!\!A_0\!\!\\ 0&\!\!A_0\!\!&0\rix.
\]
Multiplying the first and the last block-row by $-1$ we obtain a strong posH linearization
\begin{equation}\label{lin1}
L(\lambda)=\lambda (J_1+R_1)+J_2+R_2=
\lambda \mat{ccc}0&\!\!-A_3\!\!&0\\ \!A_3\!&A_2&0\\0&0&A_0\rix\!+\mat{ccc}A_3&0&0\\ 0&A_1&\!\!A_0\!\!\\ 0&\!\!-A_0\!\!&0\rix.
\end{equation}
with the coefficient matrices
\begin{equation}\label{jjrr}
J_1=\mat{ccc}0&\!\!-A_3\!\!&0\\ \!A_3\!&0&0\\0&0&0\rix\!,\;R_1=\mat{ccc}0&0&0\\ 0&\!A_2\!&0\\ 0&0&\!\!A_0\!\!\rix\!,\;
J_2=\mat{ccc}0&0&0\\ 0&0&\!\!A_0\!\!\\ 0&\!\!-A_0\!\!&0\rix\!,\;R_2=\mat{ccc}\!A_3\!&0&0\\ 0&\!A_1\!&0\\ 0&0&0\rix\!.
\end{equation}
The pencil $\lambda J_1+J_2$ is singular, since for any $\lambda_0\in\mathbb C$, the matrix $\lambda_0 J_1+J_2$ is
rank deficient, but $J_1$ and $J_2$ do not have a common left or right kernel, so all left and right minimal indices are
larger than zero. 
This means that we will not be able to apply Theorem~\ref{thm:lefthalfplane2}. However, Theorem~\ref{thm:lefthalfplane} will
lead to the following sufficient condition for the spectrum to be in the closed left half plane.
%
%
\begin{theorem}\label{Benner-explained}
Let $L(\lambda)=\lambda^3 A_3+\lambda^2 A_2+\lambda A_1 +A_0$ be a complex matrix polynomial with $A_3,A_2,A_0>0$, $A_1\geq0$
and $A_2+A_1>0$. Then the spectrum lies outside the set
\[
\big\{z\in\Comp\,\big|\, \RE z>0,\ -\beta_* < \IM z < \beta_* ,\ \arctan(-\beta_*) < \arg z< \arctan(\beta_*)\big\},
\]
where
\[
\beta_*=\sup\left\{\beta\geq 0\,\left|\; \mat{cc} A_3 & -\beta \ii A_3\\ \ii\beta A_3 & A_1+A_2 \rix>0\right\}\right..
\]
If, additionally,
%
%
%
\begin{equation}\label{pos2}
A_2\geq A_3\ \mbox{\rm  and}\ A_1\geq A_0,
\end{equation}
then all eigenvalues of $L(\lambda)$ lie in the closed left half-plane.
\end{theorem}
\proof
We use the strong posH linearization from~\eqref{lin1} and show that this pencil has the desired properties.
The first part of the proof then follows directly from Theorem~\ref{3/4}\eqref{WCiii}, note that due to the form of the
matrices in \eqref{jjrr} we have $\beta_*=\beta_+=\beta_-$.

Furthermore, with the matrices in~\eqref{jjrr} in the pencil $L(\lambda)$ in \eqref{skewlin} we have that $R_1,R_2\geq 0$ and
$\ker R_1\cap\ker R_2=\set0$. Hence, $J_1+R_1$ and $J_2+R_2$ do not have a common isotropic vector.
In order to apply Theorem~\ref{thm:lefthalfplane}, it remains to show that
that the condition \eqref{EE-JJ} is satisfied. Setting $x=\mat{ccc} x_1^T & x_2^T & x_3^T \rix^T \in\Comp^{3n}$,
we can write \eqref{EE-JJ} as
\begin{equation}\label{AA-DD}
-(x_2^*A_2 x_2+x_3^*A_0x_3)(x_1A_3x_1^* +x_2^* A_1x_2) +  4 \RE (x_2^* A_3 x_1)\RE( x_2^* A_0 x_3)  \leq 0.
\end{equation}
%
It remains to show that \eqref{pos2} implies \eqref{AA-DD}.
Setting $a^{(k)}_i=x_i^* A_kx_i$ for $i,k=1,2,3$, then  \eqref{pos2} implies that $a_i^{(2)}\geq a_i^{(3)}, a_i^{(1)}\geq a_i^{(0)}$
for $i=1,2,3$.  Hence
\[
a_1^{(3)}(a_2^{(2)}-a_2^{(3)})+a^{(0)}_3(a^{(1)}_2-a^{(0)}_2) +a^{(2)}_2 a^{(1)}_2 - a^{(0)}_2a^{(3)}_2\geq 0,
\]
and we obtain that
\begin{eqnarray*}
(a^{(2)}_2+a^{(0)}_3)(a^{(3)}_1+a^{(1)}_2)  \geq & (a^{(0)}_2+a^{(3)}_1)(a^{(3)}_2+a^{(0)}_3)\\
\geq & 4 ( a^{(3)}_1 a^{(3)}_2 a^{(0)}_2 a^{(0)}_3)^{1/2}\\
\geq & 4|x_2^* A_3 x_1| |x_2^* A_0 x_3|,
\end{eqnarray*}
which shows \eqref{AA-DD}, where in the second inequality we used the inequality for geometric and arithmetic means and in the
third inequality we used the Cauchy-Schwarz inequality.
\eproof
\begin{example} \rm {A simple calculation shows that the scalar polynomial $\lambda^3+a\lambda^2+a\lambda+1$ has roots in the
closed left half-plane if and only if $a\geq 1$. Hence, condition 
\eqref{pos2} is sharp.}
\end{example}

The following corollary applies to pencils appearing in Moore-Gibson-Thompson eigenvalue equation, see \cite{Ben21,KalN19}.
\begin{corollary}\label{Benner-explained2}
Let $L(\lambda)=I_n\lambda^3+a I_n\lambda^2+ bT\lambda + cT$, where $a,b,c>0$ and $T>0$.
If $a>1$ and $b>c$ then all eigenvalues lie in the left half-plane.
\end{corollary}

For a matrix polynomial in \eqref{PSC} of arbitrary degree $d\geq 3$  it seems  difficult to obtain general statements under which
conditions on the coefficients of matrix polynomial all eigenvalues are contained in the closed left half plane.
Hence one may have to check for each case individually which results from Section~\ref{sec3} can be applied to a posH linearization
of the given matrix polynomial.

\section*{Conclusion}
We have studied (posH) matrix pencils with coefficients having positive semidefinite Hermitian parts and  matrix polynomials
with Hermitian positive semidefinite coefficients. These generalize dissipative Hamiltonian pencils or are their matrix polynomial analogues.
We have characterized when posH pencils are equivalent to dissipative Hamiltonian pencils and we have presented several results
that lead to restrictions for the spectral properties of posH pencils and  matrix polynomials with Hermitian positive semidefinite
coefficients. This includes, in particular, the singular part and the parts associated with infinite eigenvalues.

\section*{Acknowledgment}
The authors are indebted to \L ukasz Kosi\'nski for interesting discussions on complex analysis, which inspired some  results from Section~\ref{sec4}.

\bibliographystyle{plain}

\end{document}